\documentclass[12pt]{article}
\usepackage{amssymb}
\usepackage{amscd}
\usepackage{amsmath}
\newtheorem{thm}{Theorem}[section]

\newtheorem{lem}[thm]{Lemma}
\newtheorem{cor}[thm]{Corollary}
\newtheorem{defn}[thm]{Definition}

\newtheorem{theorem}[thm]{Theorem}
\newtheorem{proposition}[thm]{Proposition}
\newtheorem{lemma}[thm]{Lemma}
\newtheorem{corollary}[thm]{Corollary}
\newtheorem{definition}[thm]{Definition}
\newtheorem{remark}[thm]{Remark}

\def\Td{\mathop{\operatorname{\rm Td}\nolimits}}
\def\Id{\mathop{\operatorname{\rm Id}\nolimits}}

\def\Ind{\mathop{\operatorname{\rm Ind}\nolimits}}
\def\Mat{\mathop{\operatorname{\rm Mat}\nolimits}}
\def\GL{\mathop{\operatorname{\rm GL}\nolimits}}

\def\HP{\mathop{\operatorname{\rm HP}\nolimits}}
\def\Hom{\mathop{\operatorname{\rm Hom}\nolimits}}
\def\Tot{\mathop{\operatorname{\rm Tot}\nolimits}}
\def\ad{\mathop{\operatorname{\rm ad}\nolimits}}
\def\Ad{\mathop{\operatorname{\rm Ad}\nolimits}}
\def\Aff{\mathop{\operatorname{\rm Aff}\nolimits}}
\def\aff{\mathop{\operatorname{\rm aff}\nolimits}}
\def\Aut{\mathop{\operatorname{\rm Aut}\nolimits}}
\def\Lie{\mathop{\operatorname{\rm Lie}\nolimits}}

\def\Ad{\mathop{\mathrm {Ad}}\nolimits}
\def\ad{\mathop{\mathrm {ad}}\nolimits}
\def\Lie{\mathop{\operatorname {Lie}}\nolimits}
\def\Aff{\mathop{\mathrm {Aff}}\nolimits}
\def\aff{\mathop{\mathrm {aff}}\nolimits}

\def\Id{\mathop{\mathrm {Id}}\nolimits}
\def\Ad{\mathop{\mathrm {Ad}}\nolimits}
\def\ad{\mathop{\mathrm {ad}}\nolimits}
\def\Ind{\mathop{\mathrm {Ind}}\nolimits}
\def\Id{\mathop{\mathrm {Id}}\nolimits}
\def\Aut{\mathop{\mathrm {Aut}}\nolimits}

\def\sgn{\mathop{\mathrm {sgn}}\nolimits}
\def\Lie{\mathop{\mathrm {Lie}}\nolimits}
\def\Ln{\mathop{\mathrm {Ln}}\nolimits}

\def\Aff{\mathop{\mathrm {Aff}}\nolimits}
\def\aff{\mathop{\mathrm {aff}}\nolimits}

\def\Id{\mathop{\mathrm {Id}}\nolimits}
\def\Ad{\mathop{\mathrm {Ad}}\nolimits}
\def\ad{\mathop{\mathrm {ad}}\nolimits}
\def\Ind{\mathop{\mathrm {Ind}}\nolimits}
\def\Id{\mathop{\mathrm {Id}}\nolimits}
\def\Aut{\mathop{\mathrm {Aut}}\nolimits}

\def\sgn{\mathop{\mathrm {sgn}}\nolimits}

\large
\begin{document}
\title{A Survey of Noncommutative Chern Characters}
\author{Do Ngoc Diep\footnote{The work was supported in part by the Alexander von Humboldt Foundation, Germany, the Abdus Salam International Center for Theoretical Physics (ICTP), Italy and National Foundation for Research in Fundamental Sciences, Vietnam, and it was completed during a visit of the author at The University of Iowa, U.S.A.}}
\maketitle
\begin{abstract}
We report in this survey some new results concerning noncommutative Chern characters: construction and the cases when they are exactly computed. The major result indicates some clear relation of these noncommutative objects and their commutative counterparts. This survey can be considered as the second part of the previous survey \cite{diep1}.
\end{abstract}
\tableofcontents
\section{Introduction}

\subsection{Motivation}

It is well-known the Chern characters $$ch : K^*(X) \to H^*_{DR}(X)$$ as homomorphisms from cohomological K-groups to $\mathbf Z/(2)$-graded de Rham cohomology groups. The first theory is related with index theory and the second theory - with integration over manifolds.
It is also well-known that modulo torsions, they are isomorphisms. More precisely, without torsions, they become isomorphisms $$ch_\mathbf{Q} : K^*(X)\otimes \mathbf{Q} \to H^*_{DR}(X;\mathbf{Q}).$$ This means that we can consider the homology theory as some approximation of K-theory. And therefore, the indices of elliptic operators related with integration over manifolds. The discovery of Atiyah-Hirzebruch-Singer index formulla reflects this relation, expressing the indices of elliptic operators, which count the differences of dimensions of kernels and cokernels of Fredholm operators, with some integrals over manifolds, which evaluate some characteristic classes over the fundamental class of manifolds. 

Chern characters play an important role in the problem of describing the structure of group C*-algebras and quantum groups, particularly in the Atiyah-Hirzebruch-Singer index formula for Dirac operators on $Spin^c$-manifold,
$$\Ind \mathbf D = 2^{\frac{d(d-1)}{2}}\langle \Td(X) \cup Ch_E, [X]\rangle, \mbox{ where } d = \dim(X). $$

It is well-known that the algebras of functions on manifolds define the structure of manifolds. One can therefore define manifolds as some special class of commutative algebras. From this we have therefore a new look to the notion of manifolds: In the K-theory context this can be expressed as follows. Following a well-known Serre-Swan theorem, the topological K-theory of the topological compact $X$ is isomorphic to the algebraic K-theory $K_*(\mathbf{C}(X))$ of the Banach algebra $\mathbf{C}(X)$ of continuous fuctions on $X$. In the de Rham theory context, the cohomology of topological spaces is isomorphic to the de Rham homology of the algebras of smooth functions on manifolds. The Chern characters have sense as homomorphisms, and in many cases, isomorphisms between them.
Our main goal is to do ``the same'' for more general (noncommutative) algebras. In this survey, three kinds of noncommutative algebras are considered: compact Lie group C*-algebras, compact quantum group C*-algebras and some classes of quantum algebras via deformation quantization:

The group C*-algebras were introduced many year ago, but still it is difficult to describe the structure of the group C*-algebras of noncompact locally compact groups. We introduce the new idea related with the Chern characters of noncommutative group C*-algebras.
We are going also to consider algebras of functions on the quantum groups. These algebras are obtained from some algebraic deformation as it was in the works of Vaksman and Soibelman. A rich variety of quantum algebras is the deformation quantization of Poisson structure. The idea simulizes what was done in Physics.  

\subsection{Tools and Ideas}

Noncommutative geometry with two important tools: KK-theory and cyclic homology. A. Connes introduced the cyclic homology theory and then the entire cyclic periodic homology $HE_*(A)$. He shows that this theory has all good properties of a homology theory like: homotopy invariant, Morita invariance and excision and constructed the Chern characters
$$ch : K_*(A)\to HE_*(A)$$ with the help of pairing $$K_*(A) \times HE^*(A) \to \mathbf{C}.$$ This theory recovered the classical results when $A=\mathbf{C}(X)$ is a commutative Banach algebra. It is more difficult to do any computation for noncommutative algebras. We consider the case when the algebra $A$ can be presented as an inductive limits of good ideals $I_\alpha$ and restrict the pairing of Connes
$$K_*(A) \times HE^*(A) \to \mathbf{C}$$ to every ideal $I_\alpha$ with $ad$-invariant trace
$$K_*(A) \times HE^*(I_\alpha)\to \mathbf{C}$$ and then extend to $$K_*(A) \times \cup_\alpha HE^*(I_\alpha) \to \mathbf{C}.$$ This therefore gives us a homomorphisms $$K_*(A) \to (\cup_\alpha HE^*(I_\alpha))'$$ as some generalized Chern character. The question is to introduc a modified entire periodic cyclic homology $HE_*(\varinjlim I_\alpha)$ as some new $HE_*(A)$ with the same properities.
This modification, fall constructed gives us a possibility to compute the noncommutative Chern characters in the indicated 3 cases: the group C*-algebras of compact Lie groups, the compact quantum group C*-algebras and some quantum coadjoint orbits, appeared from deformation quantization. 
This review is closely related to and can be considered as the second part of our previous review\cite{diep1} and book\cite{diep2}.

\section{Compact Lie Group C*-Algebras}
Let us recall $G$ a compact Lie group, $L^1(G)$ involutive Banach algebra of functions with convolution product w. r. t. the natural Haar measure. The Banach $L^1$ norm is irregular, i.e. in general $$\Vert a\* a^* \Vert \ne \Vert a \Vert^2 .$$ One introduced, therfore, the regular norm $$\Vert a\Vert_{C^*(G)}:= \sup_{\pi\in \hat{G}} \Vert \pi(a)\Vert.$$ The C*-algebra is the completion with respect to this C* -norm. 
It is well known for compact groups:
\begin{enumerate}
\item Every representations is unitarizable and every irreducible unitary representation is finite dimensional, namely of dimension $n_i$ 
\item The set of equivalent classes of irreducible unitary representations are not more than denombrable.
\end{enumerate}

As consequence it was proven:

\begin{thm} The group C*-algebras of compact groups can be expressed as  topological Cartisean product of matix algebras: $$C^*(G) \cong (\widetilde{\prod})_{i=1}^\infty Mat_{n_i}(\mathbf C),$$
\end{thm}

Let us recall that \cite{dt1}
$\tau_\alpha$ is an $ad$-invariant trace on an ideal $I_\alpha$ in a Banach algebra $A$ iff:
\begin{enumerate}
\item $\tau_\alpha$ is continuous linear, $\Vert\tau_\alpha\Vert =1$
\item $\tau_\alpha$ is positive $\tau_\alpha(a^*a) \geq 0, \forall \alpha$ and strictly positive $\tau_\alpha(a^*a) = 0 \mbox{ iff } a = 0$
\item $\tau_\alpha$ is $ad$-invariant, i.e. $\tau_\alpha([a,x]) = 0, \forall x\in A , a\in \tau_\alpha.$
\end{enumerate}
\begin{cor} The C*-algebra $C^*(G)$ can be presented as inductive limit of ideal
$$C^*(G) \cong \varinjlim_N I_N, \qquad I_N := \prod_{i=1}^N Mat_{n_i}(\mathbf{C}).$$
\end{cor}

\section{Algebras of Functions on Compact Quantum  Groups}
For a compact Lie group $G$ with Lie algebra $\mathfrak{g}$, define Hopf algebra $\mathcal{F}_\varepsilon(G)$ of functions over the quantum group, which is dual to the quantized universal enveloping algebra $U_\varepsilon(\mathfrak{g})$.

\begin{defn}{\rm
C*-norm:
$$\Vert f\Vert := \sup_{\rho} \Vert\rho(f) \Vert  \quad (f\in \mathcal{F}_\varepsilon(G)) $$
}\end{defn}

\begin{thm}
$$C^*_\varepsilon(G) \cong \mathbf{C}(\mathbf{T})\bigoplus_{e\ne w\in W}\int_{W\times \mathbf T} {\mathcal K}(\mathbf H_{w,t}) dt,$$
\end{thm}
\begin{cor}
Compact quantum group C*-algebras can be presented as inductive limits of ideals with $ad$-invariant trace.
\end{cor}

\section{C*-Algebraic Noncommutative Chern  Characters}
\subsection{K-Theory of Banach algebras}
$A$ Banach algebra with unit 1
	
\begin{thm} For Banach *-algebras,
$$KK(A,\mathbf{C}) \cong K_*^{top}(A)\cong K_*^{alg}(A)$$
\end{thm}
Let us recall that the algebraic K*-groups are defined as follows: $K_0(A)$ is the Abelian group envelop of the Grothendieck semi-group, generated by projective finitely generated $A$-modules. The $K_1$-groups are defined as the Whitehead groups, which are the abelized infinite general linear groups $\GL_\infty(A)/[\GL_\infty(A), \GL_\infty(A)]$

\subsection{Entire homology}

Let us consider a Banach *-algebra which can be presented as a direct limit of special ideal $I_\alpha$ 
$$A = \varinjlim_\alpha (I_\alpha, \tau_\alpha),$$
where
$\tau_\alpha$ is an $ad$-invariant trace on $I_\alpha$, i.e. :
\begin{enumerate}
\item $\tau_\alpha$ is continuous linear, $\Vert\tau_\alpha\Vert =1$
\item $\tau_\alpha$ is positive $\tau_\alpha(a^*a) \geq 0, \forall \alpha$ and strictly positive $\tau_\alpha(a^*a) = 0 \mbox{ iff } a = 0$
\item $\tau_\alpha$ is $ad$-invariant, i.e. $\tau_\alpha([a,x]) = 0, \forall x\in A , a\in \tau_\alpha.$
\end{enumerate}
We then have for every $\alpha \in \Gamma$ a scalar product $$\langle
a,b\rangle_\alpha
:=\tau_\alpha(a^*b)$$ and also an direct system $\{I_\alpha,
\tau_\alpha\}_{\alpha \in
\Gamma}$. Let $\bar{I}_\alpha$ be the completion of $I_\alpha$ with respect to the scalar product defined
above and $\widetilde{\bar{I}_\alpha}$ denote $\bar{I}_\alpha$ with formally
adjointed unity
element. Define $C^n(\widetilde{\bar{I}_\alpha})$ as the set of
$n+1$-linear maps $\varphi :
(\widetilde{\bar{I}_\alpha})^{\otimes(n+1)} \to {\mathbf C}$. There exists a Hilbert structure on  $(\widetilde{\bar{I}_\alpha})^{\otimes(n+1)}$ and we can identify
$C_n((\widetilde{\bar{I}_\alpha})) := \Hom(C^n(\bar{\tilde{I}}_\alpha), {\mathbf C})$ with $C^n(\bar{\tilde{I}}_\alpha)$ via an anti-isomorphism.

For $I_\alpha \subseteq I_\beta$, we have a
well-defined map $$D^\beta_\alpha :
C^n(\widetilde{\bar{I}_\alpha}) \to C^n(\widetilde{\bar{I}_\beta}),$$ which
makes $\{
C^n(\widetilde{\bar{I}}_\alpha)\}$ into a direct system. Write $Q = \varinjlim
C^n(\widetilde{\bar{I}_\alpha})$ and observe that it admits a Hilbert
space structure, see \cite{dt1}-\cite{dt2}. Let $C_n(A) := \Hom(\varinjlim
C^n(\widetilde{\bar{I}_\alpha}), {\mathbf C}) = \Hom(\varinjlim
C_n(\widetilde{\bar{I}_\alpha}), {\mathbf C})$ which is anti-isomorphic to
$\varinjlim_\alpha C^n(\bar{\tilde{I}}_\alpha)$. So we have finally
$$C_n(A) = \varinjlim C_n(\bar{\tilde{I}}_\alpha).$$

Let $$b, b' : C^n(\widetilde{\bar{I}}_\alpha) \to
C^{n+1}(\widetilde{\bar{I}}_\alpha),$$
$$N : C^n(\widetilde{\bar{I}}_\alpha) \to
C^n(\widetilde{\bar{I}}_\alpha),$$ $$\lambda :
C^n(\widetilde{\bar{I}}_\alpha) \to C^n(\widetilde{\bar{I}}_\alpha),$$
$$S : C^{n+1}(\widetilde{\bar{I}}_\alpha)\to
C^n(\widetilde{\bar{I}}_\alpha)$$ be defined as in A. Connes \cite{Ca}.
We adopt the notations in \cite{Kh1} . Denote by $b^*, (b')^*, N^*,
\lambda^*, S^*$ the corresponding adjoint operators. Note also that for
each
$I_\alpha$ we have the same formulae for adjoint operators for homology as
Connes obtained for cohomology.

We now have a bi-complex
$$
\leqno{{\mathcal C}(A):}  
\CD
 @.    \vdots  @.    \vdots @.     \vdots @.      \\
@.          @V(-b')^* VV        @Vb^* VV       @V(-b')^* VV        @.   \\
\dots @<1-\lambda^* << \varinjlim_\alpha C_1(\bar{\tilde{I}}_\alpha) @<N^*<<
\varinjlim_\alpha C_1(\bar{\tilde{I}}_\alpha)
@<1-\lambda^*<< \varinjlim_\alpha C_1(\bar{\tilde{I}}_\alpha) @<N^*<< \dots\\
@.          @V(-b')^* VV        @Vb^* VV       @V(-b')^* VV        @.   \\
\dots @<1-\lambda^*<< \varinjlim_\alpha C_0(\bar{\tilde{I}}_\alpha) @<N^*<<
\varinjlim_\alpha 
C_0(\bar{\tilde{I}}_\alpha) @<1-\lambda^*<<\varinjlim_\alpha
C_0(\bar{\tilde{I}}_\alpha) @<N^*<< \dots\\
\endCD$$
with $d_v := b^*$ in the even columns and $d_v:= (-b')^*$ in the odd columns, $d_h
:= 1-\lambda^*$ from odd to even columns and $d_h := N^*$ from even to odd columns, 
where * means the corresponding adjoint operator. Now we have
$$\Tot({\mathcal C}(A))^{even} = \Tot({\mathcal C}(A))^{odd} :=
\oplus_{n\geq 0}
C_n(A),$$ which is periodic with period two. Hence, we have
$$\begin{array}{ccc}
&\partial & \\
\oplus_{n\geq 0} C_n(A)&\begin{array}{c} \longleftarrow\\ \longrightarrow
\end{array} &\oplus_{n\geq 0} C_n(A)\\
&\partial & \\
\end{array}$$
where $\partial = d_v + d_h$ is the total differential.

\begin{defn} Let $\HP_*(A)$ be the homology of the total complex
$(\Tot{\mathcal C}(A))$. It is called the {\sl periodic cyclic homology} 
of $A$.
\end{defn}

Note that this $\HP_*(A)$ is, in general different from the 
$\HP_*(A)$ of Cuntz-Quillen, because we used the direct limit of periodic
cyclic homology of ideals. But in special cases, when the whole
algebra $A$ is one of these ideals
 with $\ad_A$-invariant trace, (e.g. the
commutative algebras
of complex-valued functions on compact spaces), we return to the
Cuntz-Quillen $\HP_*$, which we shall use later. 

\begin{defn}
An even (or odd) chain $(f_n)_{n\geq 0}$ in ${\mathcal C}(A)$ is called
{\it entire} if the
radius of convergence of the power series $\sum_n \frac{n!}{[\frac{n}{2}]!}
\Vert
f_n\Vert z^n$, $z\in {\mathbf C}$ is infinite.
\end{defn}

Let $C^e_*(A)$ be the sub-complex of $C(A)$ consisting of entire chains. Then
we have a
periodic complex.

\begin{thm} Let
$$\Tot(C^e_*(A))^{even} = \Tot(C^e_*(A))^{odd} := \bigoplus_{n\geq 0} C^e_n(A),$$
where $C^e_n(A)$ is the entire $n$-chain. Then we have a complex of
entire chains with the total differential $\partial := d_v + d_h$
$$\begin{array}{ccc} & \partial
& \\
\bigoplus_{n\geq 0} C^e_n(A) & \begin{array}{c}\longleftarrow \\
\longrightarrow\end{array} &
\bigoplus_{n\geq 0}
C^e_n(A)\\
& \partial & \end{array}$$
\end{thm}
\begin{defn}
The homology of this complex is called also the {\it entire homology} and
denoted by $HE_*(A)$. 
\end{defn}
Note that this entire homology is defined through
the inductive limits of ideals with ad-invariance trace. 

In \cite{dt1}-\cite{dt2}, the main properties of this theory, namely
\begin{itemize}
\item Homotopy invariance,
\item Morita invariance and
\item Excision,
\end{itemize}
were proven and hence $HE_*$ is a generalized homology theory.

\begin{lem}
If the Banach algebra $A$ can be presented as a direct limit
$\varinjlim_\alpha I_\alpha$ of a system
of ideals $I_\alpha$ with $\ad$-invariant trace $\tau_\alpha$, then
$$K_*(A) = \varinjlim_\alpha K_*(I_\alpha),$$
$$HE_*(A) = \varinjlim_\alpha HE_*(I_\alpha).$$
\end{lem}
The following result from K-theory is well-known:
\begin{thm} The entire homology of non-commutative de Rham currents
admits the following stability property
$$K_*({\mathcal K}({\mathbf H})) \cong K_*({\mathbf C}),$$
$$K_*(A \otimes {\mathcal K}({\mathbf H})) \cong K_*(A),$$
where ${\mathbf H}$ is a separable Hilbert space and $A$ is an arbitrary Banach
space.
\end{thm}

The similar result is true for entire homology $HE_*$ :
\begin{thm} The entire homology of non-commutative de Rham currents
admits the following stability property
$$HE_*({\mathcal K}({\mathbf H})) \cong HE_*({\mathbf C}),$$
$$HE_*(A \otimes {\mathcal K}({\mathbf H})) \cong HE_*(A),$$
where ${\mathbf H}$ is a separable Hilbert space and $A$ is an arbitrary Banach
space.
\end{thm}
	
Let $A$ be an involutive Banach algebra. In this section, we construct a
non-commutative
character
$$ch : K_*(A) \to HE(A)$$ and later show that when $A= C^*(G)$, this Chern
character
reduces up to isomorphism to classical Chern character.

Let $A$ be an involutive Banach algebra with unity. \begin{thm}
There exists a Chern character $$ch : K_*(A) \to HE_*(A).$$ \end{thm}

We first recall that there exists a pairing $$K_n(A) \times C^n(A) \to
{\mathbf C}$$ due to
A. Connes, see \cite{Co}. Hence there exists a map
$ K_n(A) \stackrel {C_n }{ \longrightarrow }\Hom(C^n(A),{\mathbf C}).$ So,
by 1.1,
we have for each $\alpha \in \Gamma$ a pairing
$$K_n(A) \times C^n(\tilde{\bar{I}}_\alpha) \to \mathbf{C}$$
 and hence a map $K_n(A)
\stackrel{C_n^\alpha}{\longrightarrow}
\Hom(C^n(\widetilde{\bar{I}_\alpha}),{\mathbf C})$ and hence a map $K_n(A)
\stackrel{C_n}{\longrightarrow}
\Hom(C^n(\widetilde{\bar{I}_\alpha}),{\mathbf C}).$ We now show that $C_n$
induces the Chern map
$$ch: K_n(A) \to HE_n(A)$$

Now let $e$ be an idempotent in $M_k(A)$ for some $k\in {\mathbf N}$. It
suffies to
show that for $n$ even, if $\varphi = \partial\psi$, where $\varphi\in
C^n(\widetilde{\bar{I}_\alpha})$ and $\psi \in
C^{n+1}(\widetilde{\bar{I}_\alpha})$,
then $$\langle e,\varphi\rangle = \sum_{n=1}^\infty
\frac{(-1)^n}{n!}\varphi(e,e,\dots,e)
= 0.$$
However, this follows from Connes' results in (\cite{Co}, Lemma 7).

The proof of the case for $n$ odd would also follow from \cite{Co}. 

\subsection{ Compact Lie group C*-algebra case}

Our next result computes the Chern character in 2.1 for $A=C^*(G)$ by
reducing it to the classical case.
Let us recall that the group C*-algebras of compact locally compact groups can be presented as some inductive limits of $ad$-invariant ideals of type $I_N = \prod_{i=1}^N Mat_{n_i}(\mathbf{C})$,
$$C^*(G) = \varinjlim_N I_N, \mbox{ where } I_N = \prod_{i=1}^N Mat_{n_i}(\mathbf{C}).$$

Let ${\mathbf T}$ be a fixed maximal torus of $G$ with Weyl group $W:=
\mathcal{N}_G({\mathbf T})/{\mathbf T}$.
\begin{thm}
Then the Chern character $$K^G_*({\mathbf C}^*(G)) \to HE^G_*(C^*(G))$$ is an
isomorphism, which can be identified with the classical Chern character
$$K^W_*({\mathbf C}({\mathbf T})) \to HE^W_*({\mathbf C}({\mathbf T}))$$ that
is also an isomorphism.
\end{thm}

Because of the pairing
$$K_*(A) \times HE^*(A) \to \mathbf{C},$$
and because $A_\alpha$ are supposed to be ideals in $A$, 
we can extend it to the pairing
$$K_*(A) \times \varinjlim_\alpha HE^*(A_\alpha) \to \mathbf{C},$$
and therefore
$$K_*(A) \times \bigcup_\alpha\varinjlim_\alpha HE^*(A_\alpha) \to \mathbf{C}.$$

\subsection{Quantum group C*-algebra case}

In this case of compact quantum group C*-algebras we have also
$$C^*_\varepsilon(G) \cong \mathbf{C}(\mathbf T) \oplus \int^\oplus_{\mathcal{N}_\mathbf{T}\setminus \{e\} \times \mathbf T} \mathcal{K}(\mathbf{H}_{w,t}) dt$$

Our next result computes the Chern characters  for $A=C^*_\varepsilon(G)$ by
reducing it to the classical case.

\begin{thm}
Let ${\mathbf T}$ be a fixed maximal torus of $G$ with Weyl group $W:=
{\mathcal N}_{\mathbf T}/{\mathbf T}$.
Then the Chern character $$ch_{C^*}: K_*(C^*_\varepsilon(G)) \to
HE_*(C^*_\varepsilon(G))$$ is an
isomorphism modulo torsion, i.e.
$$\CD
ch_{C^*}: K_* (C^*_\varepsilon(G))\otimes {\mathbf C} @>\cong >>
HE_*(C^*_\varepsilon(G)),\endCD $$
which can be identified with the classical
Chern character
$$\CD ch: K_*({\mathbf C}({\mathcal N}_{\mathbf T})) @>>> HE_*({\mathbf
C}({\mathcal N}_{\mathbf T}))\endCD$$
that
is also an isomorphism modulo torsion, i.e.
$$\CD
ch: K_* ({\mathcal N}_{\mathbf T})\otimes {\mathbf C} @>\cong >>
H^*_{DR}({\mathcal N}_{\mathbf T}).
\endCD$$ 
\end{thm}

\section{Quantum Strata of Coadjoint Orbits}

Let us now consider the third class of noncommutative algebras, arized from deformation quantization. In this section we expose the results obtained in \cite{diep4}.
For locally compact groups, their C*-algebras contain exhausted informations about the groups them-selves and their representations, see \cite{diep1}, \cite{diep2}. In some sense \cite{rieffel1}-\cite{rieffel2}, the group algebras can be considered as C*-algebraic deformation quantization $C^*_q(G)$ at the special value $q=1$.

In \cite{diep1} and \cite{diep2}, the group C*-algebras were described as repeated extensions of C*-algebras of strata of coadjoint orbits. Quantum groups are group Hopf algebras, i.e. replace C*-algebras by Hopf algebras. It is therefore interesting to ask whether we could describe quantum groups as some repeated extensions of some kind quantum strata of coadjoint orbits? We are attempting to give a positive answer to this question. It is not yet completely described but we obtained a reasonable answer. Let us describe the main ingredients of our approach.

For the good strata of coadjoint orbits, there exist always continuous fields of polarizations (in the sense of the representation theory), satisfying the L. Pukanszky irreducibility condition: for each orbit ${\mathcal O}$ and any point $F_{\mathcal O}$ in it, the affine subspace, orthogonal to polarizations with respect to the symplectic form  is included in orbits themselves,  i.e.
$$F_{\mathcal O} + \mathfrak h_\mathcal{O}^\perp \subseteq \mathcal{O}$$
and $$\dim \mathfrak{h}_{\mathcal O} = \frac{1}{2}\dim \mathcal{O}.$$ We choose the the canonical Darboux coordinates with impulse $p$'s-coordinates, following a vector structure basis of $\mathfrak h^\perp$. From this we can deduce that in this kind of Darboux coordinates, the Kirillov form $\omega_{\mathcal O}$ locally are canonical and every element $X\in \mathfrak g = \Lie G$ can be considered as a function $\tilde{X}$ on $\mathcal O$, linear on $p$'s-coordinates, i.e. $$\tilde{X} = \sum_{i=1}^n a_i(q)p_i + a_0(q).$$ 

This essential fact gives us a possibility to effectively write out the corresponding $\star$-product of functions, and define quantum strata $C_q(V,\mathcal F)$. On the strata acts our Lie group of symmetry. It induces therefore an action on equivariant differential operators. Using the indicate fields of polarizations we prove some kind of Poincar\'e-Birkhoff-Witt theorem and then provide quantization with separation of variable in sense of Karabegov \cite{karabegov1}. We can then express the corresponding representations of the quantum strata $C_q(V,\mathcal F)$, where $V \subseteq \bigcup_{\dim \mathcal O = const} \mathcal O$, through the Feynman path integrals etc...,see \cite{diep3}. 

\subsection{Canonical coordinates on a stratum}

Let us consider a connected and simply connected Lie group $G$ with Lie algebra $\mathfrak g$. Denote the dual to $\mathfrak g$ vector space by $\mathfrak g^*$. It is well-known that the action of $G$ on itself by conjugation $$A(g) : G \to G,$$ defined by $A(g)(h) := ghg^{-1}$ keeps the identity element $h=e$ unmoved. This induces the tangent map $\Ad(g) := A(g)_*: {\mathfrak g} = T_eG\to \mathfrak g,$ defined by $$\Ad(g)X := \frac{d}{dt}|_{t=0} A(g)\exp(tX)$$
and the co-adjoint action 
$K(g) := A(g^{-1})^* $ 
maps the dual space 
${\mathfrak g}^*$ 
into itself. The orbit space 
${\mathcal O}(G):= {\mathfrak g}^*/G$ 
is in general a bad topological space, namely non-Hausdorff. Consider 
an arbitrary orbit 
$\Omega\in{\mathcal O}(G)$
and an element 
$F \in {\mathfrak g}^*$ in it. The stabilizer is denote by 
$G_F$, 
its connected component by $(G_F)_0$ and its Lie algebra by 
${\mathfrak g}_F := \Lie(G_F)$. It is well-known that
$$
\begin{array}{ccc}
G_F & \hookrightarrow & G\\
    &                 & \Big\downarrow\\
    &                  & \Omega_F
    \end{array}
$$
is a principal bundle with the structural group $G_F$. 
Let us fix some {\it connection in this principal bundle, } 
\index{connection on 
principal bundle} i.e. some {\it trivialization } \index{trivialization}
of this bundle.  
We want to construct representations in some cohomology spaces 
with coefficients in the sheaf of sections of some vector bundle 
associated with this principal bundle. 
It is well know  that every vector bundle 
is an induced one with respect to some representation of the structural group
in the typical fiber.
It is natural to fix some unitary representation 
$\tilde{\sigma}$ 
of $G_F$ such that its kernel contains $(G_F)_0$, the character
$\chi_F$ of the connected component of stabilizer
$$\chi_F(\exp{X}) := \exp{(2\pi\sqrt{-1}\langle F,X\rangle )}$$
and therefore the differential 
$D(\tilde{\sigma}\chi_F)=\tilde{\rho}$ 
is some representation  of the Lie algebra ${\mathfrak g}_F$. 
We suppose that the representation $D(\tilde{\rho}\chi_F)$ was extended to 
the complexification $({\mathfrak g}_F)_{\mathbf C}$. 
The whole space of all sections seems to be so large for the construction 
of irreducible unitary representations. 
One consider the invariant subspaces with the help of some so 
called polarizations, see \cite{diep1}, \cite{diep2}.

\begin{definition} We say that a triple $(\mathfrak p_\mathcal O, \rho_\mathcal O, \sigma_{0,\mathcal O})$ is a 
$(\tilde{\sigma},F)$-polarization of $\mathcal O$ iff :
\begin{enumerate}
\item[a.] $\mathfrak p_\mathcal O$ is some complex sub-algebra of the complexified $\mathfrak g_\mathbf C$, containing $\mathfrak g_{F_\mathcal O}$.
\item[b.] The sub-algebra $\mathfrak p_\mathcal O$ is invariant under the action of all the operators of type $Ad_{\mathfrak g_\mathbf C} x, x\in G_{F_\mathcal O}.$
\item[c.] The vector space $\mathfrak p_\mathcal O + \overline{\mathfrak p_\mathcal O}$ is the complexification of some real Lie sub-algebra $\mathfrak m_\mathcal O := (\mathfrak p_\mathcal O + \overline{\mathfrak p}_\mathcal O) \cap \mathfrak g.$
\item[d.] All the subgroups $M_{0, \mathcal O}$, $H_{0, \mathcal O}$, $M_\mathcal O$, $H_\mathcal O$ are closed, where by definition, $M_{0,\mathcal O}$ (resp., $H_{0, \mathcal O}$) is the connected subgroup of $G$ with the Lie algebra $\mathfrak m_\mathcal O$ (resp., $\mathfrak h_\mathcal O := \mathfrak p_\mathcal O \cap \mathfrak g$) and the semi-direct products $M:= M_{0,\mathcal O} \ltimes G_{F_\mathcal O}$, $H_\mathcal O := H_{0,\mathcal O} \ltimes G_{F_\mathcal O}$. 
\item[e.] $\sigma_{0,\mathcal O}$ is an irreducible representation of $H_{0,\mathcal O}$ in some Hilbert space $V_\mathcal O$ such that: 
(1.) the restriction $\sigma_{0,\mathcal O}\vert_{G_{F_\mathcal O} \cap H_{0,\mathcal O}}$ is some multiple of the restriction $\chi_{F_\mathcal O}.\tilde{\sigma}_\mathcal O\vert_{G_{F_\mathcal O} \cap H_{0,\mathcal O}}$, where the character $\chi_\mathcal O$ is by definition, $\chi_\mathcal O(\exp X) = \exp(2\pi\sqrt{-1}\langle F_\mathcal O, X\rangle )$;
(2.) under the action of $G_{F_\mathcal O}$ on the dual $\hat{H}_{0,\mathcal O}$, the point $\sigma_{0,\mathcal O}$ is fixed.
\item[f.] $\rho_\mathcal O$ is a representation of the complex Lie algebra $\mathfrak p_\mathcal O$ in the same $V_\mathcal O$, which satisfies the Nelson conditions for $H_{0,\mathcal O}$ and $\rho_\mathcal O \vert_{\mathfrak h_\mathcal O} \cong D\sigma_{0,\mathcal O}$.
\end{enumerate}
\end{definition}
There is a natural order in the set of all $(\tilde{\sigma}_\mathcal O, F_\mathcal O)$-polarizations by inclusion and from now on speaking about polarizations we means always the maximal ones. It is not hard to prove that the (maximal) polarizations are also the Lagrangian distributions and in particular the co-dimension of $\mathfrak h_\mathcal O$ in $\mathfrak g$ is a half of the dimension of the coadjoint orbit $\mathcal O$,  $$codim_\mathfrak g \mathfrak h_\mathcal O = \frac{1}{2}\dim \mathcal O,$$ see e.g. \cite{diep2} or \cite{kirillov1}.

Let us now recall the Pukanszky condition.
\begin{definition}
We say that the $(\tilde{\sigma}_\mathcal O,F_\mathcal O)$-polarization $(\mathfrak p_\mathcal O, \rho_\mathcal O, \sigma_{0,\mathcal O})$ satisfies the Pukanszky condition, iff $$F_\mathcal O + \mathfrak h_\mathcal O^\perp \subset \mathcal O.$$
\end{definition}

\begin{remark}
\begin{itemize}
\item Pukanszky conditions involve an inclusion of the Lagrangian affine subspace of $p$'s coordinates into the local Darboux coordinates.
\item The partial complex structure on orbits let us to use smaller subspaces of section in the induction construction, as subspaces of partially invariant, partially holomorphic sections of the induced bundles.
\end{itemize}
\end{remark}

\begin{thm}
There exists on each coadjoint orbit a local canonical system of Darboux coordinates, in which the Hamiltonian function $\tilde{X}$, $X \in \mathfrak g$, are linear on $p$'s impulsion coordinates and in theses coordinates, $$\tilde{X} = \sum_{i=1}^n a_i(q)p_i + a_0(q).$$
\end{thm}

\subsection{Poisson structure on strata}

Let us first recall the construction of strata of coadjoint orbit from
\cite{diep1}-\cite{diep2}. The orbit space $\mathcal O(G)$ is a disjoint union of $\Omega_{2n}$, each of which is a union of the coadjoint orbits of dimension $2n$, $ 0 \leq 2n \leq \dim G$. Denote $$V_{2n} = \cup_{\dim \mathcal O = 2n} \mathcal O.$$ Then $V_{2n}$ is the set of points of fixed rank $2n$ of the Poisson structure bilinear function $$\{X,Y\} = \langle F,[X,Y]\rangle .$$

Suppose that it is a foliation, at least for $V_{2n}$, where $2n$ is the maximal dimension possible in $\mathcal O(G)$. It can be shown that the foliation $V_{2n}$ can be obtained from the group action of $\mathbf R^{2n}$ on $V_{2n}$. Let us for this aim, consider a basis $X_1, \dots, X_{2n}$ of the tangent space $T_{F_\mathcal O}\mathcal O \cong \mathfrak g /\mathfrak g_{F_\mathcal O}$ at the point $F_\mathcal O \in \mathcal O \subset V_{2n}$. We can define an action of $\mathbf R^{2n}$ on $V_{2n}$ as 
$$(t_1,\dots,t_{2n}) \mapsto \exp(t_1X_1)\exp(t_2X_2)\dots \exp(t_{2n}X_{2n})F_\mathcal O .$$
We have therefore the Hamiltonian vector fields
$$\xi_k := \frac{d}{dt}\vert_{t=0} \exp(t_kX_k)F, \forall k=1,\dots, 2n$$ and their span $\mathcal F_{2n} = \{\xi_1,\dots,\xi_{2n}\}$ provides a tangent distribution. It is easy to show that we have therefore a measurable (in sense of A. Connes) foliation. One can therefore define also the Connes C*-algebra $C^*(V_{2n},\mathcal F_{2n})$, $0 \leq 2n \leq \dim G$. By introducing some technical condition in \cite{diep1}, it is easy to reduced these C*-algebras to extensions of other ones, those are in form of tensor product $C(X) \otimes \mathcal K(H)$ of algebras of continuous functions on compacts  and the elementary algebra $\mathcal K(H)$ of compact operators in a separable Hilbert space $H$. The strata of these kinds we means $good$ strata. Another kind of strata of coadjoint orbits are obtained from relation with the cases where the Gelfand-Kirillov conjecture was solve, for examples for connected and simply connected solvable Lie groups, see \S5.

It is deduced from a result of Kontsevich that this Poisson structure can be quantized. This quantization however is formal. The question of convergence of the quantizing series is not clear. We show in this section the case of charts with the linear p's canonical coordinates, the corresponding $\star$-product  is convergent.

In this kind of special local chart systems of Darboux coordinates it is easy to deduce existence local convergent $\star$-products. 
\begin{thm}
Locally on each coadjoint orbit, there exist a $star$-product.
\end{thm}
 
Let us denote by $\mathcal F^{-1}_p$ the Fourier inverse transformation on variables $p$'s and by $\mathcal F_p$ the Fourier transformation. For two symbols $f,g \in \mathcal F^{-1}_p(PDO_G(\mathcal O))= \mathbf C_q(\mathcal O)$, their Fourier images $\mathcal F_p(f)$ and $\mathcal F_p(g)$ are operators and we define their oprator product and then take Fourier inverse, as $\star$-product
$$f \star g := \mathcal F^{-1}_p(\mathcal F_p(f).\mathcal F_p(g)).$$ So this product is again a symbol in the same class $\mathbf C_q(\mathcal O)$.

Because of existence of special coordinate systems, linear on $p$'s coordinates we can treat for the good strata in the same way as in the cases of exponential groups.  And the formal series of $\star$-product is convergent.
The proof could be also done in the same scheme as in exponential or compact cases, see \cite{arnalcortet1}-\cite{arnalcortet2}.

Let us denote by $\Gamma= \pi_1({\mathcal O})$ the fundamental group of the orbit. Our next step is to globally extend this kind of local $\star$-products. Our idea is as follows. We lift this $\star$-product to the universal covering of coadjoint orbits as some $\Gamma$-invariant $\star$-products, globally extend them in virtue of the monodromy theorem and then pushdown to our orbits.
We start with the following lemmas

\begin{lem} \label{lemma1}
There is one-to-one correspondence between  $\star$-products on Poisson manifolds and $\Gamma$-invariant $\star$-products on their the universal coverings.
\end{lem}

We use this lemma to describe existence of a $\star$-product on coadjoint orbits.

\begin{lem}\label{lemma2}
On a universal covering, each local $\star$-product can be uniquely extended to some global $\star$-product on this covering. 
\end{lem}
 For local charts, there exist deformation quantization, as said above, $f \mapsto Op(f)$ by using the formulas of Fedosov quantization. Also in the intersection of two local charts of coordinates $(q,p)$ and $(\tilde{q}, \tilde{p})$, there is a symplectomorphism, namely $\varphi$ such that $$(\tilde{q},\tilde{p}) = \varphi(q,p).$$ Using the local oscilatting integrals and by compensating the local Maslov's index obstacles, one can exactly construct the unitary operator $U$ such that $$Op(f\circ \varphi) = U.Op(f) .U^{-1},$$ see for example Fedosov's book \cite{fedosov}.
Because the universal coverings are simply connected, the extensions can be therefore produced because of Monodromy Theorem.

\begin{thm}\label{mainthm}
There exists a convergent $\star$-product on each orbit, which is a symplectic leaf of the Poisson structure on each stratum of coadjoint orbits.
\end{thm}

From the description of the canonical coordinates in the previous section, we see that there exist locally $\star$-product. This local $\star$-product then extended to  a $\Gamma$-invariant global one on the universal covering, which produces a convergent $\star$-product on the coadjoint orbits, following Lemmas \ref{lemma1},\ref{lemma2}.

\subsection{Star-product, Quantization and PBW Theorem}

We use the constructed in the previous section $\star$-product to provide an action of functions on the spaces of partially invariant partially holomorphic sections of the corresponding partially invariant holomorphically induced bundles associated with polarizations. It is possibles because the quantum induced bundles are locally trivial and the spaces of partially invariant partially holomorphic sections with section in Hilbert spaces are locally finitely generated modules over the algebras of quantizing functions, see \cite{diep2}. 
We use then the construction of Karabegov and Fedosov to obtain a Hopf *-algebra of longitudinal pseudo-differential operators elliptic along the leaves of this measurable foliations. The main ingredient is that we use here the Poincar\'e-Birkhopf-Witt theorem to provide this quantization. 

\begin{thm}\label{theorem3.1}
There is a natural deformation quantization with separation of variables, corresponding to the Poincar\'e-Birkhoff-Witt Theorem for polarizations.
\end{thm}

We have from the multidimensional quantization $\mathfrak g \to PDO^1(\mathcal O)$, $X \mapsto \hat{X}$. More precisely, Let us denote by $PDO(\mathcal O)$ the algebra of right $G$-invariant pseudo-differential operators on $\mathcal O$, i.e. the continuous maps from $C^\infty(\mathcal O)$ into itself not extending support. We also denote $PDO^1(\mathcal O)$ the Lie algebras of right $G$-invariant first order pseudo-differential operators. By the procedure of multidimensional quantization, there is a homomorphism of Lie algebras
$$\mathfrak g \to PDO_G^1(\mathcal O) \subset PDO_G(\mathcal O).$$ Following the universal property of $U(\mathfrak g)$, there is a unique homomorphism of associative algebras $U(\mathfrak g) \to PDO_G(\mathcal O)$making the following diagram commutative
$$\CD
\mathfrak g @>>> PDO_G(\mathcal O)\\
@VVV    @AAA\\
U(\mathfrak g) @= U(\mathfrak g)
\endCD$$

\begin{remark}
Because of BKW, $U(\mathfrak g) \cong U(\mathfrak p_{\mathcal O}/\mathfrak h_{\mathcal O}) \otimes U(\mathfrak h_{\mathcal O}) \otimes U(\overline{\mathfrak p}_{\mathcal O}/\mathfrak h_{\mathcal O})$ and because of this our quantizing map is coincided with that one used by Karabegov in the quantization with separation of variables in case of coadjoint orbits with totally complex polarizations $\mathfrak g_\mathbf C = \mathfrak m_\mathbf C = \mathfrak p \oplus \overline{\mathfrak p}$.
\end{remark}

 Let us first describe the machinery applied in the method of Karabegov's separation of variable. We use the idea about polarizations in multidimensional quantization, \cite{diep2}.

Let us recall the root decomposition $$\mathfrak g = \mathfrak n_- \oplus \mathfrak a \oplus \mathfrak n_+.$$ If $\mathfrak p$ is a complex polarization, then from definition we have
$\mathfrak m_\mathbf C = \mathfrak p \oplus \bar{\mathfrak p},$
$\mathfrak h_\mathbf C = \mathfrak p \cap \bar{\mathfrak p},$
where $\mathfrak m = \mathfrak g \cap (\mathfrak p \oplus \bar{\mathfrak p}) $ and 
$\mathfrak h = \mathfrak g \cap (\mathfrak p \oplus \bar{\mathfrak p}).$ The quotients subspaces $\mathfrak p / \mathfrak h_{\mathbf C}$ and $\bar{\mathfrak p}/\mathfrak h_{\mathbf C}$ are included in $\mathfrak m$ as linear subspaces (not necessarily to be sub-algebras). Let us fix some (non-canonical) inclusions.
Let us denote by $U(\mathfrak p/\mathfrak h_{\mathbf C})$ (resp. $U(\bar{\mathfrak p}/\mathfrak h_{\mathbf C})$) the {\it sub-algebra, generated} by elements from $\mathfrak p/\mathfrak h_{\mathbf C} \hookrightarrow \mathfrak m_{\mathbf C}$ (resp. $\mathfrak p/\mathfrak h_{\mathbf C} \hookrightarrow \mathfrak m_{\mathbf C}$) in the universal algebra $U(\mathfrak m)$. We have therefore a analog of the well-known Poincar\'e-Birkhoff-Witt theorem.
\begin{thm}[Poincar\'e-Birkhoff-Witt Theorem]
If  $\mathfrak p$ is as polarization, then $$U(\mathfrak m_{\mathbf C}) \cong U(\mathfrak p/\mathfrak h_{\mathbf C}) \otimes U(\mathfrak h_{\mathbf C}) \otimes U(\bar{\mathfrak p} /\mathfrak h_{\mathbf C})$$ 
\end{thm}

Let us fix in a basis each of $\mathfrak p/\mathfrak h_\mathbf C$, $\mathfrak h_\mathbf C$ and $\bar{\mathfrak p}/\mathfrak h_\mathbf C$. We have therefore a basis of $\mathfrak m_\mathbf C = \mathfrak p /\mathfrak h_\mathbf C \oplus \mathfrak h_\mathbf C \oplus \bar{\mathfrak p}/\mathfrak h_\mathbf C$. Our theorem is therefore deduced from the original Poincar\'e-Birkhoff-Witt Theorem \cite{pbw}.

It is easy to see that by this reason, on the variety $M/H$ there is a natural complex structure and therefore on coadjoint orbits exists some partial complex structure. We use this complex structure and this Poincar\'e-Birkhoff-Witt theorem to do a separation of variables on $M$ and apply the machinery of Karabegov.
\begin{thm} \label{separationthm}
In the case of totally complex polarizable coadjoint orbits, the quantization map from the theorem \ref{theorem3.1} is coincided with the rule of Karabegov's quantization with separation of variables.
\end{thm}

\subsection{Representations}

Let us consider now a continuous fields of $(\tilde{\sigma}_\mathcal O,F_\mathcal O)$-polarizations $(\mathfrak p_\mathcal O, \rho_\mathcal O, \sigma_{0,\mathcal O})$ satisfying the Pukanszky condition. On one hand side, we can use the multidimensional quantization procedure to obtain irreducible unitary representations $\Pi_\mathcal O =  Ind(G;\mathfrak p_\mathcal O,H_\mathcal O,\rho_\mathcal O,\sigma_{0,\mathcal O})$ of $G$, \cite{diep2}. On other hand side, we can use $\star$-product construction to provide the representations $\Pi_\mathcal O : \tilde{X\vert_\Omega} \mapsto \ell_X$ of the quantum strata $\mathbf C_q(\Omega)$. We'd like to show we have the same one.

\begin{defn}
{\it Quantum coadjoint orbit} $\mathbf C_q(\mathcal O)$ is defined as the Hopf algebra of symbols of differential operators $U_{q,\mathcal O}(\mathfrak g) = U(\mathfrak g)\vert_{\mathcal O}$. The homomorphism $Q: U(\mathfrak g) \to PDO_G(\mathcal O)$ is defined to be the {\it second quantization homomorphism}.
\end{defn}

\begin{thm}
 The representation of the Lie algebra obtained from $\star$-product is equal to the representations obtained from the multidimensional quantization procedure.
\end{thm}

 Let us recall \cite{diep2}, that 
$$\Lie_X \; Ind(G;\mathfrak p_\mathcal O,H_\mathcal O,\rho_\mathcal O,\sigma_{0,\mathcal O}) \cong \hat{X}.$$
From the construction of quantization map $$U(\mathfrak g) \to PDO_G(\mathcal O)$$ as the map arising from the universal property the map $\mathfrak g \to PDO^1_G(\mathcal O),$ $$\ell_X = Op(\tilde{X}) = \hat{X}, \forall X \in \mathfrak g = \Lie(G).$$ The associated representation of $\mathbf C_q(\Omega)$ is the solution of the Cauchy problem for the differential equation
$$\frac{\partial}{\partial t}U(t, q,p) = \ell_X U(t,q,p),$$
$$U(0,q,p) = \Id .$$
The solution of this problem is uniquely defined.

\subsection{Oscillating Fourier integrals}

Let us in this section consider the good family of coadjoint orbits arising from the solved cases of Gelfand-Kirillov Conjecture. The results in this section is revised from the \cite{diep3}.

Consider a connected and simply connected Lie group $G$ with Lie algebra $\mathfrak g$.
\begin{thm}
There exists a $G$-invariant Zariski open set $\Omega$ and a covering $\tilde{\Omega}$ of $\Omega$, with natural action of $G$ such that for each continuous field of polarizations $(\mathfrak p_\mathcal O, H_\mathcal O, \rho_\mathcal O, \sigma_{0,\mathcal O})$, $\mathcal O \in \Omega/G$, the Lie derivative of the direct integral of representations arized from the multidimensional quantization procedure is equivalent to the tensor product of the Schr\"odinger representation $\Pi = Sch$ of the Gelfand-Kirillov basis of the enveloping field and a continuous field of trivial representations $\{V_\mathcal O\}$ on $\tilde{\Omega}$.
\end{thm}

The proof is rather long and consists of several steps:

1. We apply the construction of Nghiem \cite{nghiem5} to the solvable radical $^r{\mathfrak g}$ of $\mathfrak g$ to obtain the Zariski open set $^r\Omega$ in $^r\mathfrak g^*$ and its covering $^r\tilde{\Omega}$.

2. The general case is reduced to the semi-simple case $^s\mathfrak g$, see (\cite{nghiem1}, Thm. C). Denote by $^sG$ the corresponding analytic subgroup of $G$.

3. Take the Zariski open set $\mathcal A_s\mathcal P(^sG)$ of admissible and well-polarizable strongly regular functionals from $^s\mathfrak g^*$ and its covering $\mathcal B_s(^sG)$ via Duflo's construction \cite{duflo1}.

4. The desired $G$-invariant Zariski open set and its covering are the corresponding Cartesian products
$$\Omega= \mathcal A_s\mathcal P(^sG) \times {}^r\Omega^0 \times {}^r\Omega^1 \times \dots \times {}^r\Omega^k$$
$$\tilde{\Omega}= \mathcal B_s(^sG) \times {}^r\tilde{\Omega^0} \times {}^r\tilde{\Omega^1} \times \dots\times  {}^r\tilde{\Omega}^k$$

5. \begin{lem}
There exists a continuous field of polarizations of type $(\mathfrak p_\mathcal O, H_\mathcal O,\rho_\mathcal O,\sigma_{0,\mathcal O})$, for each $\mathcal O \in \tilde{\Omega}/G$.
\end{lem}

6. \begin{lem}
The Lie derivative commutes with direct integrals, i.e. $$\nabla \int^\oplus_{\tilde{\Omega}/G}Ind(G;\mathfrak p_\mathcal O, H_\mathcal O,\rho_\mathcal O,\sigma_{0,\mathcal O})d\mathcal O = \int^\oplus_{\tilde{\Omega}/G} \nabla\;Ind(G;\mathfrak p_\mathcal O, H_\mathcal O,\rho_\mathcal O, \sigma_{0,\mathcal O}) d\mathcal O.$$
\end{lem}

7. \begin{lem}
The restriction of the Schr\"odinger representation to coadjoint orbits provides a continuous field of polarizations. In particular,
$$D\sigma_{0,\mathcal O} = \rho_\mathcal O\vert_{\mathfrak p_{\mathcal O} \cap \mathfrak g} \cong mult \; Sch_\mathcal O.$$
\end{lem}

\begin{thm}\label{fourierthm}
1. There exists an operator-valued phase $\Phi(t,z,x)$ and an operator-valued amplitude $a(t,z,x,y,\xi)$ extended from the expression
$$\exp(\Phi(t,z,y) + \sqrt{-1}\xi(g(t)x-x))$$ in such a fashion that the action of the representation $$\Pi_\mathcal O = Ind(G;\mathfrak p_\mathcal O, H_\mathcal O,\rho_\mathcal O,\sigma_{0,\mathcal O})$$ can be expressed as an oscilatting Fourier integral
$$\Pi_\mathcal O(g(t))f(z,x) = c.\int_{\mathbf R^M}\int_{\mathbf R^M} a(t,z,x,y,\xi)\exp(\sqrt{-1}\xi(x-y))f(z,y)dyd\xi,$$ where $c$ is some constant.

2. For each function $\varphi$ of Schwartz class $\mathcal S(G)$ satisfying the compactness criteria \cite{diep2} in every induction step \cite{lipsman}, \cite{diep2}, the operator $\Pi_\mathcal O(\varphi)$ is of trace class and its action can be expressed as the oscilatting Fourier integral 
$$\Pi_\mathcal O(\varphi)f(z,x) = const.\int_{\mathbf R^M}\int_{\mathbf R^M}(\int_{\mathbf R^\ell} a(t,z,x,y,\xi)\varphi(g(t))dt) \times$$ $$\times \exp(\sqrt{-1}\xi(x-y))f(z,y)dyd\xi. $$
Hence, its trace is $$tr\Pi_\mathcal O(\varphi) = const.\int_{\mathbf R^M}\int_{\mathbf R^M}(\int_{\mathbf R^\ell} a(t,z,x,x,\xi)\varphi(g(t))dt )dxd\xi$$
\end{thm}

The proof also consists of several steps:

1. From \cite{diep1} - \cite{diep2} and \cite{duflo1} it is easy to see obtain a slight unipotent modification( i.e. a reduction to the unipotent radical) of the multidimensional quantization procedure. We refer the reader to \cite{tranvui1}-\cite{tranvui2}, and \cite{trandong} for a detailed exposition.

2. From the unitary representations $\Pi_\mathcal O = Ind(G; \mathfrak p_\mathcal O, H_\mathcal O, \rho_\mathcal O , \sigma_{0,\mathcal O})$ in the unipotent context and its differential $\pi_\mathcal O = D\Pi_\mathcal O(G;  \mathfrak p_\mathcal O, H_\mathcal O, \rho_\mathcal O, \sigma_{0,\mathcal O})$, it is easy to select the constant term and the vector fields term, (see \cite{nghiem4} for the simplest case and notations), 
$$\pi(L) = \pi_c(L) + \pi_d(L), \quad \pi_s(L) = a_0(L, x),$$
$$\pi_d(L) = \partial (L) = \sum_k a_k(L) \frac{\partial}{\partial x_k} .$$

The operator-valued partial (i.e. depending on $L$) phase $$\Phi_L(t_L,z,x) = \int_0^{t_L} a_0(z,g(s)x)ds.$$ The phase $\Phi(t,z,x)$ is the sum of the partial phases, on which the one-parameter group $g_L(s)$ operates by translations for all the factors on the left of $g_L(t_L)$ in the ordered product $$g(t) = \prod_L g_L(t_L). $$

It is easy then to see that our induced representations $\Pi_\mathcal O$ act as follows
$$\Pi_\mathcal O(g(t))f(z,x) = \sigma_{0,\mathcal O}(t,z,x) f(z,g(t)x).$$

3. It suffices now to apply the Fourier transforms
$$f(z,x) = (2\pi )^M \int_{\mathbf R^M} \int_{\mathbf R^M} \exp\{\sqrt{-1}\xi (x-y)\} f(z,y)dyd\xi$$ to the function $\Pi_\mathcal O(g(t))f(z,x)$ to obtain
$$\Pi_\mathcal O(g(t))f(z,x) = $$
$$(2\pi)^M\int_{\mathbf R^M}\int_{\mathbf R^M} a(t,z,x,y,\xi)\exp\{\sqrt{-1}\xi (x-y)\}f(z,y)dyd\xi,$$
where the amplitude $a(t,z,x,y,\xi)$ is the natural extension of the expression $$\exp\{ \Phi(t,z,y) + \sqrt{-1}\xi (g(t)x -x) \}$$ in the correspondence with the fields of polarizations. Hence for each $\varphi \in \mathcal S(G)$,
$$\Pi_\mathcal O(\varphi)f(z,x) = (2\pi)^M \int_{\mathbf R^M} \int_{\mathbf R^M} (\int_{\mathbf R^\ell} a(t,z,x,y\xi)\varphi(g(t))dt) \times$$
$$\times \exp\{\sqrt{-1}\xi (x - y) \}f(z,y)dyd\xi .$$

Remark that the integral $\int_{\mathbf R^\ell} a(t,z, x,y,\xi)\varphi(g(t)) dt$ is just a type of Feynman path integrals.

4. \begin{lem}
If in every repeated induction step, see \cite{lipsman}, \cite{diep2}, $\Pi_\mathcal O(\varphi)$ satisfies the compactness criteria, then the  operator $\Pi_\mathcal O(\varphi)$ is trace class and hence $$tr\; \Pi_\mathcal O(\varphi) = (2\pi)^M \int_{\mathbf R^M}\int_{\mathbf R^M} tr( \int_{\mathbf R^\ell} a(t,z,x,x,\xi)\varphi(g(t))dtdxd\xi.$$
\end{lem}

\section{Examples}
In this section, we illustrate the main ideas in examples.
\subsection{$\overline{MD}$-groups}
We expose in this subsection our joint works with Nguyen Viet Hai \cite{diephai1}-\cite{diephai2}.
\subsubsection{The group of affine transformations of the real straight line}
We refer the reader to the work \cite{diephai1} for a detailed exposition with complete proof.

{\it Canonical coordinates on the upper half-planes.}
Recall that the Lie algebra $\mathfrak g = \aff(\mathbf  R)$ of affine transformations of the real straight line is described as follows, see for example \cite{diep2}: The Lie group $\Aff(\mathbf R)$ of affine transformations of type $$x \in \mathbf R \mapsto ax + b, \mbox{ for some parameters }a, b \in \mathbf R, a \ne 0.$$ It is well-known that this group $\Aff(\mathbf R)$ is a two dimensional Lie group which is isomorphic to the group of matrices
$$\Aff(\mathbf R) \cong \{\left (\begin{array}{cc} a & b \\ 0 & 1 \end{array} \right) \vert a,b \in \mathbf R , a \ne 0 \}.$$ We consider its connected component $$G= \Aff_0(\mathbf R)= \{\left (\begin{array}{cc} a & b \\ 0 & 1 \end{array} \right) \vert a,b \in \mathbf R, a > 0 \}$$ of identity element. Its Lie algebra is
$$\mathfrak g = \aff(\mathbf R) \cong  \{\left (\begin{array}{cc} \alpha & \beta \\ 0 & 0 \end{array} \right) \vert \alpha, \beta  \in \mathbf R \}$$  admits a basis of two generators $X, Y$ with the only nonzero Lie bracket $[X,Y] = Y$, i.e. 
$$\mathfrak g = \aff(\mathbf R) \cong \{ \alpha X + \beta Y \vert [X,Y] = Y, \alpha, \beta \in \mathbf R \}.$$
The co-adjoint action of $G$ on $\mathfrak g^*$ is given (see e.g. \cite{arnalcortet2}, \cite{kirillov1}) by $$\langle K(g)F, Z \rangle = \langle F, \Ad(g^{-1})Z \rangle, \forall F \in \mathfrak g^*, g \in G \mbox{ and } Z \in \mathfrak g.$$ Denote the co-adjoint orbit of $G$ in $\mathfrak g$, passing through $F$ by 
$$\Omega_F = K(G)F :=  \{K(g)F \vert F \in G \}.$$ Because the group $G = \Aff_0(\mathbf R)$ is exponential (see \cite{diep2}), for $F \in \mathfrak g^* = \aff(\mathbf R)^*$, we have 
$$\Omega_F = \{ K(\exp(U)F | U \in \aff(\mathbf R) \}.$$
It is easy to see that
$$\langle K(\exp U)F, Z \rangle = \langle F, \exp(-\ad_U)Z \rangle.$$ It is easy therefore to see that
$$K(\exp U)F = \langle F, \exp(-\ad_U)X\rangle X^*+\langle F, \exp(-\ad_U)Y\rangle Y^*.$$
For a general element $U = \alpha X + \beta Y \in \mathfrak g$, we have
$$\exp(-\ad_U) = \sum_{n=0}^\infty \frac{1}{n!}\left(\begin{array}{cc}0 & 0 \\ \beta & -\alpha \end{array}\right)^n = \left( \begin{array}{cc} 1 & 0 \\ L & e^{-\alpha} \end{array} \right),$$ where $L = \alpha + \beta + \frac{\alpha}{\beta}(1-e^\beta)$. This means that
$$K(\exp U)F = (\lambda + \mu L) X^* + (\mu e^{\-\alpha})Y^*. $$ From this formula one deduces  \cite{diep2} the following description of all co-adjoint orbits of $G$ in $\mathfrak g^*$:
\begin{itemize}
\item If $\mu = 0$, each point $(x=\lambda , y =0)$ on the abscissa ordinate corresponds to a 0-dimensional co-adjoint orbit $$\Omega_\lambda = \{\lambda X^* \}, \quad \lambda \in \mathbf R .$$
\item For $\mu \ne 0$, there are two 2-dimensional co-adjoint orbits: the upper half-plane $\{(\lambda , \mu) \quad\vert\quad \lambda ,\mu\in \mathbf R , \mu > 0 \}$ corresponds to the co-adjoint orbit
\begin{equation} \Omega_{+} := \{ F = (\lambda + \mu L)X^* + (\mu e^{-\alpha})Y^* \quad \vert \quad \mu > 0 \}, \end{equation}
and the lower half-plane $\{(\lambda , \mu) \quad\vert\quad \lambda ,\mu\in \mathbf R , \mu < 0\}$ corresponds to the co-adjoint orbit
\begin{equation} \Omega_{-} := \{ F = (\lambda + \mu L)X^* + (\mu e^{-\alpha})Y^* \quad \vert \quad \mu < 0 \}. \end{equation}
\end{itemize}
We shall work from now on for the fixed co-adjoint orbit $\Omega_+$. The case of the co-adjoint orbit $\Omega_-$ is similarly treated. First we study the geometry of this orbit and introduce some canonical coordinates in it.
It is well-known from the orbit method \cite{kirillov1} that the Lie algebra $\mathfrak g = \aff(\mathbf R)$, realized by the complete right-invariant Hamiltonian vector fields on co-adjoint orbits $\Omega_F \cong G_F \setminus G$ with flat (co-adjoint) action of the Lie group $G = \Aff_0(\mathbf R)$. On the orbit $\Omega_+$ we choose a fix point $F=Y^*$. It is well-known from the orbit method that we can choose an arbitrary point $F$ on $\Omega_F$. It is easy to see that the stabilizer of this (and therefore of any) point  is trivial $G_F = \{e\}$. We identify therefore $G$ with $G_{Y^*}\setminus G$. There is a natural diffeomorphism $\Id_{\mathbf R} \times \exp(.)$ from the standard symplectic space $\mathbf R^2$ with symplectic 2-form $dp \wedge dq$ in canonical Darboux $(p,q)$-coordinates, onto the upper half-plane $\mathbf H_+ \cong \mathbf R \rtimes
 \mathbf R_+$ with coordinates $(p, e^q)$, which is, from the above coordinate description, also diffeomorphic to the co-adjoint orbit $\Omega_+$. We can use therefore $(p,q)$ as the standard canonical Darboux coordinates in $\Omega_{Y^*}$. There are also non-canonical Darboux coordinates $(x,y) = (p,e^q)$ on $\Omega_{Y^*}$. We show now that in these coordinates $(x,y)$, the Kirillov form looks like $\omega_{Y^*}(x,y) = \frac{1}{y}dx \wedge dy$, but in the canonical Darboux coordinates $(p,q)$, the Kirillov form is just the standard symplectic form $dp \wedge dq$. This means that there are  symplectomorphisms between the standard symplectic space $\mathbf R^2, dp \wedge dq)$, the upper half-plane $(\mathbf H_+, \frac{1}{y}dx \wedge dy)$ and the co-adjoint orbit $(\Omega_{Y^*},\omega_{Y^*})$. 
Each element $Z\in \mathfrak g$ can be considered as a linear functional $\tilde{Z}$ on co-adjoint orbits, as subsets of $\mathfrak g^*$, $\tilde{Z}(F) :=\langle F,Z\rangle$.  It is well-known that this linear function is just the Hamiltonian function associated with the Hamiltonian vector field $\xi_Z$, which represents $Z\in \mathfrak g$ following the formula 
$$(\xi_Zf)(x) := \frac{d}{dt}f(x\exp (tZ))|_{t=0}, \forall f \in C^\infty(\Omega_+).$$ 
The Kirillov form $\omega_F$ is defined by the formula 
\begin{equation}\label{7} \omega_F(\xi_Z,\xi_T) = \langle F,[Z,T]\rangle, \forall Z,T \in \mathfrak g = \aff(\mathbf R). \end{equation} This form defines the symplectic structure and the Poisson brackets on the co-adjoint orbit $\Omega_+$. For the derivative along the direction $\xi_Z$ and the Poisson bracket we have relation $\xi_Z(f) = \{\tilde{Z},f\}, \forall f \in C^\infty(\Omega_+)$. It is well-known in differential geometry that the correspondence
$Z \mapsto \xi_Z, Z \in \mathfrak g$ defines a representation of our Lie algebra by vector fields on co-adjoint orbits. If the action of $G$ on $\Omega_+$ is flat \cite{diep2}, we have the second Lie algebra homomorphism from  strictly Hamiltonian right-invariant vector fields into the Lie algebra of smooth functions on the orbit with respect to the associated Poisson brackets.

Denote by $\psi$ the indicated symplectomorphism from $\mathbf R^2$ onto $\Omega_+$
$$(p,q) \in \mathbf R^2 \mapsto \psi(p,q):= (p,e^q) \in \Omega_+$$
\begin{proposition}
1. Hamiltonian function $f_Z = \tilde{Z}$ in canonical coordinates $(p,q)$ of the orbit $\Omega_+$ is of the form $$\tilde{Z}\circ\psi(p,q) = \alpha p + \beta e^q, \mbox{ if  } Z = \left(\begin{array}{cc} \alpha & \beta \\ 0 & 0 \end{array} \right).$$

2. In the canonical coordinates $(p,q)$ of the orbit $\Omega_+$, the Kirillov form $\omega_{Y^*}$ is just the standard form $\omega = dp \wedge dq$.
\end{proposition}

{\it Computation of generators $\hat{\ell}_Z$}
Let us denote by $\Lambda$ the 2-tensor associated with the Kirillov standard form $\omega = dp \wedge dq$ in canonical Darboux coordinates. We use also the multi-index notation. Let us consider the well-known Moyal $\star$-product of two smooth functions $u,v \in C^\infty(\mathbf R^2)$, defined by
$$u \star v = u.v + \sum_{r \geq 1} \frac{1}{r!}(\frac{1}{2i})^r P^r(u,v),$$ where
$$P^r(u,v) := \Lambda^{i_1j_1}\Lambda^{i_2j_2}\dots \Lambda^{i_rj_r}\partial_{i_1i_2\dots i_r} u \partial_{j_1j_2\dots j_r}v,$$ with $$\partial_{i_1i_2\dots i_r} := \frac{\partial^r}{\partial x^{i_1}\dots \partial x^{i_r}}, x:= (p,q) = (p_1,\dots,p_n,q^1,\dots,q^n)$$ as multi-index notation. It is well-known that this series converges in the Schwartz distribution spaces $\mathcal S (\mathbf R^n)$. We apply this to the special case $n=1$. In our case we have only $x = (x^1,x^2) = (p,q)$. 
\begin{proposition}\label{3.1}
In the above mentioned canonical Darboux coordinates $(p,q)$ on the orbit $\Omega_+$, the Moyl $\star$-product satisfies the relation
$$i\tilde{Z} \star i\tilde{T} - i\tilde{T} \star i\tilde{Z} = i\widetilde{[Z,T]}, \forall Z, T \in \aff(\mathbf R).$$
\end{proposition}

Consequently, to each adapted chart $\psi$ in the sense of \cite{arnalcortet2}, we associate a $G$-covariant $\star$-product.

\begin{proposition}[see \cite{gutt}]
Let $\star$ be a formal differentiable $\star$-product on $C^\infty(M, \mathbf R)$, which is covariant under $G$. Then there exists a representation $\tau$ of $G$ in $\Aut N[[\nu]]$ such that 
$$\tau(g)(u \star v) = \tau(g)u \star \tau(g)v.$$
\end{proposition}

Let us denote by $\mathcal F_pu$ the partial Fourier transform of the function $u$ from the variable $p$ to the variable $x$, i.e.
$$\mathcal F_p(u)(x,q) := \frac{1}{\sqrt{2\pi}}\int_{\mathbf R} e^{-ipx} u(p,q)dp.$$ Let us denote by $ \mathcal F^{-1}_p(u) (x,q)$ the inverse Fourier transform. 
\begin{lemma}\label{lem3.1}
1. $\partial_p \mathcal F^{-1}_p(p.u) = i \mathcal F^{-1}_p(x.u)  $ ,

2. $ \mathcal F_p(v) = i \partial_x\mathcal F_p(v)  $ ,

3. $P^k(\tilde{Z},\mathcal F^{-1}_p(u)) = (-1)^k \beta e^q \frac{\partial^k\mathcal F^{-1}_p(u)}{\partial^kp}, \mbox{ with } k \geq 2.$
\end{lemma}

For each $Z \in \aff(\mathbf R)$, the corresponding Hamiltonian function is $\tilde{Z} =  \alpha p + \beta e^q $ and we can consider the operator $\ell_Z$ acting on dense subspace $L^2(\mathbf R^2, \frac{dpdq}{2\pi})^\infty$ of smooth functions by left $\star$-multiplication by $i \tilde{Z}$, i.e. $\ell_Z(u) = i\tilde{Z} \star u$. It is then continuated to the whole space $L^2(\mathbf R^2, \frac{dpdq}{2\pi})$. It is easy to see that, because of the relation in Proposition (\ref{3.1}), the correspondence $Z \in \aff(\mathbf R) \mapsto \ell_Z = i\tilde{Z} \star .$ is a representation of the Lie algebra $\aff(\mathbf R)$ on the space $N[[\frac{i}{2}]]$ of formal power series in the parameter $\nu = \frac{i}{2}$ with coefficients in $N = C^\infty(M,\mathbf R)$, see e.g. \cite{gutt} for more detail.

We study now the convergence of the formal power series. In order to do this, we look at the $\star$-product of $i\tilde{Z}$ as the $\star$-product of symbols and define the differential operators corresponding to $i\tilde{Z}$. It is easy to see that the resulting correspondence is a representation of $\mathfrak g $ by pseudo-differential operators. 

\begin{proposition}
For each $Z \in \aff(\mathbf R)$ and for each compactly supported $C^\infty$ function $u \in C^\infty_0(\mathbf R^2)$, we have 
$$\hat{\ell}_Z(u) := \mathcal F_p \circ \ell_Z \circ \mathcal F^{-1}_p(u) = \alpha (\frac{1}{2}\partial_q - \partial_x)u + i\beta e^{q -\frac{x}{2}}u.$$
\end{proposition}

\begin{remark}{\rm
Setting new variables $s = q - \frac{x}{2}$, $t = q + \frac{x}{2}$, we have 
\begin{equation}
\hat{\ell}_Z(u) = \alpha\frac{\partial u}{\partial s} + i\beta e^s u,
\end{equation}
e.i. $$\hat{\ell}_Z = \alpha\frac{\partial }{\partial s} + i\beta e^s ,$$ which provides a representation of the Lie algebra $\aff (\mathbf R)$. 
}\end{remark}

{\it The associate irreducible unitary representations}

Our aim in this section is to exponentiate the obtained representation $\hat{\ell}_Z$ of the Lie algebra $\aff(\mathbf R)$ to the corresponding representation of the Lie group $\Aff_0(\mathbf R)$. We shall prove that the result is exactly the irreducible unitary representation $T_{\Omega_+}$ obtained from the orbit method or Mackey small subgroup method applied to this group $\Aff(\mathbf R)$.
Let us recall first the well-known list of all the irreducible unitary representations of the group of affine transformation of the real straight line.
\begin{theorem} [\cite{gelfandnaimark}]\label{4.1}
Every irreducible unitary representation of the group $\Aff(\mathbf R)$ of all the affine transformations of the real straight line, up to unitary equivalence, is equivalent to one of the pairwise nonequivalent representations:
\begin{itemize} 
\item the infinite dimensional representation $S$, realized in the space $L^2(\mathbf R^*, \frac{dy}{\vert y\vert})$, where $\mathbf R^* = \mathbf R \setminus \{0\}$ and is defined by the formula
$$(S(g)f)(y) := e^{iby}f(ay), \mbox{ where } g = \left(\begin{array}{cc} a & b\\ 0 & 1 \end{array}\right),$$
\item the representation $U^\varepsilon_\lambda$, where $\varepsilon = 0,1$, $\lambda \in \mathbf R$, realized in the 1-dimensional Hilbert space $\mathbf C^1$ and is given by the formula
$$U^\varepsilon_\lambda(g) = \vert a \vert^{i\lambda}(\sgn a)^\varepsilon .$$
\end{itemize}
\end{theorem}
Let us consider now the connected component $G= \Aff_0(\mathbf R)$. The irreducible unitary representations can be obtained easily from the orbit method machinery.
\begin{theorem}
The representation $\exp(\hat{\ell}_Z)$ of the group $G=\Aff_0(\mathbf R)$ is exactly the irreducible unitary representation $T_{\Omega_+}$ of $G=\Aff_0(\mathbf R)$ associated following the orbit method construction, to the orbit $\Omega_+$, which is the upper half-plane $\mathbf H \cong \mathbf R \rtimes \mathbf R^*$, i. e.
+$$(\exp(\hat{\ell}_Z)f)(y) = (T_{\Omega_+}(g)f)(y) = e^{iby}f(ay),\forall f\in L^2(\mathbf R^*, \frac{dy}{\vert y\vert}), $$ where  $g = \exp Z = \left(\begin{array}{cc} a & b \\ 0 & 1 \end{array}\right).$
\end{theorem}

By analogy, we have also
\begin{theorem}
The representation $\exp(\hat{\ell}_Z)$ of the group $G=\Aff_0(\mathbf R)$ is exactly the irreducible unitary representation $T_{\Omega_-}$ of $G=\Aff_0(\mathbf R)$ associated following the orbit method construction, to the orbit $\Omega_-$, which is the lower half-plane $\mathbf H \cong \mathbf R \rtimes \mathbf R^*$, i. e.
$$(\exp(\hat{\ell}_Z)f)(y) = (T_{\Omega_-}(g)f)(y) = e^{iby}f(ay),\forall f\in L^2(\mathbf R^*, \frac{dy}{\vert y\vert}), $$ where  $g = \exp Z = \left(\begin{array}{cc} a & b \\ 0 & 1 \end{array}\right).$
\end{theorem}
\begin{remark}{\rm
1. We have demonstrated how all the irreducible unitary representation of the connected group of affine transformations could be obtained from deformation quantization. It is reasonable to refer to the algebras of functions on co-adjoint orbits with this $\star$-product as {\it quantum ones}.

2. In a forthcoming work, we shall do the same calculation for the group of affine transformations of the complex straight line $\mathbf C$. This achieves the description of {\it quantum $\overline{MD}$ co-adjoint orbits}, see \cite{diep2} for definition of $\overline{MD}$ Lie algebras.
}\end{remark}

\subsubsection{The group of affine transformations of the complex straight line}
Recall that the Lie algebra ${\mathfrak g} = \aff({\bf C})$ of affine transformations of the complex straight line is described as follows, see [D].

It is well-known that the group $\Aff({\bf C})$ is a four (real) dimensional Lie group which is isomorphism to the group of matrices:
$$\Aff({\bf C}) \cong  \left\{\left(\begin{array}{cc} a & b \\ 0 & 1 \end{array}\right) \vert a,b \in {\bf C}, a \ne 0 \right\}$$

The most easy method is to consider $X$,$Y$ as complex generators,
$X=X_1+iX_2$ and $Y=Y_1+iY_2$. Then from the relation $[X,Y]=Y$, we get$ [X_1,Y_1]-[X_2,Y_2]+i([X_1Y_2]+[X_2,Y_1]) = Y_1+iY_2$. 
This mean that the Lie algebra $\aff({\bf C})$ is a real 4-dimensional Lie algebra, having 4 generators with the only nonzero Lie brackets: $[X_1,Y_1] - [X_2,Y_2]=Y_1$; $[X_2,Y_1] + [X_1,Y_2] = Y_2$ and we can choose another basic noted again by the same letters to have more clear Lie brackets of this Lie algebra:
$$[X_1,Y_1] = Y_1; [X_1,Y_2] = Y_2; [X_2,Y_1] = Y_2; [X_2,Y_2] = -Y_1$$

\begin{remark}{\rm
The exponential map $$\exp: {\bf C}  \longrightarrow  {\bf  C}^{*} := {\bf C} \backslash \{0\}$$  giving by $z \mapsto e^z$ is just the covering map and therefore the universal covering of ${\bf} C^*$ is $\widetilde {\bf C}^* \cong {\bf C}$. As a consequence one deduces that $$\widetilde {\Aff}({\bf C}) \cong {\bf C} \ltimes{\bf C} \cong  \{(z,w) \vert z,w \in {\bf C} \}$$  with the following multiplication law:
$$(z,w)(z^{'},w^{'}) := (z+z',w+e^{z}w')$$
}\end{remark}

\begin{remark} {\rm 
The co-adjoint orbit of $\widetilde\Aff({\bf C})$ in ${\mathfrak g}^*$  passing through $F \in {\mathfrak g}^*$ is denoted by 
$$\Omega_{F} := K(\widetilde {\Aff}({\bf C})) F = \{K(g)F \vert g \in \widetilde \Aff({\bf C})\}$$
Then, (see [D]):
\begin{enumerate}
\item Each point $(\alpha,0,0,\delta)$ is 0-dimensional co-adjoint orbit $\Omega_{(\alpha,0,0,\delta)}$
\item The open set $\beta^{2}+\gamma^{2} \ne $ 0 is the single 4-dimensional co-adjoint orbit $\Omega_{F} = \Omega_{\beta^{2}+\gamma^{2} \ne 0} $. We shall also use $\Omega_{F}$ in form $\Omega_{F} \cong {\bf C} \times {{\bf C}}^*$.
\end{enumerate}
}\end{remark}

\begin{remark}{\rm 
Let us denote:
$$\mathbf H_{k} = \{w=q_{1}+iq_{2} \in {\bf C} \vert -\infty< q_1<+\infty ; 2k\pi < q_2< 2k\pi+2\pi\}; k=0,\pm1,\dots$$
$$L=\{{\rho}e^{i\varphi} \in {\bf C} \vert 0< \rho < +\infty; \varphi = 0\} \mbox{ and } {\bf C} _{k } = {\bf C} \backslash L$$
is a univalent sheet of the Riemann surface of the complex variable multi-valued analytic function $\Ln(w)$, ($k=0,\pm 1,\dots$)
Then there is a natural diffeomorphism $w \in \mathbf H_{k} \longmapsto e^{w} \in {\bf C}_k$ with each $k=0,\pm1,\dots.$ Now consider the map:
$${\bf C} \times {\bf C} \longrightarrow \Omega_F = {\bf C} \times {\bf C}^*$$
$$(z,w) \longmapsto (z,e^w),$$
with a fixed $k \in \mathbf Z$. We have a local diffeomorphism 
$$\varphi_k: {\bf C} \times {\bf H}_k \longrightarrow {\bf C} \times {\bf C}_k$$
 $$(z,w) \longmapsto (z,e^w) $$
This diffeomorphism $\varphi_k$ will be needed in the all sequel.
}\end{remark}

On ${\bf C}\times {\bf H}_k$ we have the natural symplectic form 
\begin{equation}\omega = \frac{1 }{2}[dz \wedge dw+d\overline {z} \wedge d\overline {w}],\end{equation} induced from $\mathbf C^2$.
Put $z=p_1+ip_2,w=q_1+iq_2$ and $(x^1,x^2,x^3,x^4)=(p_1,q_1,p_2,q_2) \in {\bf R}^4$, then
$$\omega = dp_1 \wedge dq_1-dp_2 \wedge dq_2.$$ The corresponding symplectic matrix of $\omega$ is 
$$ \wedge = \left(\begin{array}{cccc} 0 & -1 & 0 & 0 \\
		           1 & 0 & 0 & 0 \\
                                                0 & 0 & 0 & 1 \\
                                                0 & 0 & -1 & 0 \end{array}\right)
\mbox{   and   }          
 \wedge^{-1} = \left(\begin{array}{cccc} 0 & 1 & 0& 0 \\
		           -1 & 0 & 0 & 0 \\
                                                0 & 0 & 0 & -1 \\
                                                0 & 0 & 1 & 0 \end{array}\right)$$

We have therefore the Poisson brackets of functions as follows. With each $f,g \in {\bf C}^{\infty}(\Omega)$ 
$$\{f,g\} = \wedge^{ij}\frac{\partial f }{\partial x^i}\frac{\partial g}{\partial x^j} =  
\wedge^{12}\frac{\partial f }{\partial p_1}\frac{\partial g}{\partial q_1}+
\wedge^{21}\frac{\partial f }{\partial q_1}\frac{\partial g}{\partial p_1} +
\wedge^{34}\frac{\partial f}{\partial p_2}\frac{\partial g}{\partial q_2} +
\wedge^{43}\frac{\partial f}{\partial q_2}\frac{\partial g}{\partial p_2} = $$
$$\ \ \ \ \ \ \ =\frac{\partial f }{\partial p_1}\frac{\partial g}{\partial q_1} -
\frac{\partial f}{\partial q_1}\frac{\partial g}{\partial p_1} -
\frac{\partial f}{\partial p_2}\frac{\partial g }{\partial q_2} +
\frac{\partial f}{\partial q_2}\frac{\partial g }{\partial p_2} = $$
$$\ \ =2\Bigl[\frac{\partial f}{\partial z}\frac{\partial g}{\partial w} -
\frac{\partial f}{\partial w}\frac{\partial g}{\partial z} +
\frac{\partial f}{\partial \overline z}\frac{\partial g}{\partial \overline{w}} -
\frac{\partial f}{\partial \overline w}\frac{\partial g}{\partial \overline z}\Bigr]$$  

\begin{proposition}
Fixing the local  diffeomorphism $\varphi_k (k \in {\bf Z})$, we have: 
\begin{enumerate}
\item For any element $A \in \aff(\mathbf C)$, the corresponding Hamiltonian function $\widetilde{A}$  in local coordinates $(z,w)$ of the orbit $\Omega_F$  is of the form
$$\widetilde A\circ\varphi_k(z,w) = \frac{1}{2} [\alpha z +\beta e^w + \overline{\alpha} \overline{z} + \overline{\beta}e^{\overline {w}}]$$
\item In local coordinates $(z,w)$ of the orbit $\Omega_F$, the symplectic Kirillov form $\omega_F$ is just the standard form (1).
\end{enumerate}
\end{proposition}

{\it Computation of Operators $\hat{\ell}_A^{(k)}$.}

\begin{proposition}\label{Proposition 3.1}
With $A,B \in \aff({\bf C})$, the Moyal $\star$-product satisfies the relation:
\begin{equation} i \widetilde{A} \star i \widetilde{B} - i \widetilde{B} \star i \widetilde{A} = i[\widetilde{A,B} ]\end{equation}
\end{proposition}

For each $A \in \hbox{aff}({\bf C}$), the corresponding Hamiltonian function is 
$$\widetilde{A} = \frac{1}{2} [\alpha  z + \beta e^{w} + \overline{\alpha}\overline z + \overline{\beta} e^ {\overline w}] $$ 
and we can consider the operator  ${\ell}_A^{(k)}$  acting on dense subspace
$L^2({\bf R}^2\times ({\bf R}^2)^*,\frac{dp_1dq_1dp_2dq_2}{(2\pi)^2} )^{\infty}$
 of smooth functions by left $\star$-multiplication by $i \widetilde{A}$, i.e:
${\ell}_A^{(k)} (f) = i \widetilde{A} \star f$. Because of the relation in Proposition 3.1, we have
\begin{corollary}\label{Consequence 3.2}
 \begin{equation}{\ell}_{[A,B]}^{(k)} = {\ell}_A^{(k)} \star {\ell}_B^{(k)} - {\ell}_B^{(k)} \star {\ell}_A^{(k)} := {\Bigl[ {\ell}_A^{(k)},  {\ell}_B^{(k)}\Bigr]}^{\star}\end{equation}
\end{corollary}

From this it is easy to see that, the correspondence $A \in \aff({\bf C}) \longmapsto {\ell}_A^{(k)} = $i$\widetilde {A} \star$. is a representation of the Lie algebra $\aff({\bf C}$) on the space N$\bigl[[\frac{i}{2}]\bigr] $ of formal power series, see [G] for more detail.

Now, let us denote  ${\mathcal F}_z$(f) the partial Fourier transform of the function f from the variable $z=p_1+ip_2$ to the variable $\xi=\xi_1+i\xi_2$, i.e:
$${\mathcal F}_z(f )(\xi,w) = \frac{1 }{2\pi} \iint_{R^2} e^{-iRe(\xi \overline{z})} f(z,w)dp_1dp_2$$

Let us denote by $$\mathcal F_z^{-1}(f )(\xi,w) = \frac{1}{2\pi} \iint_{R^2} e^{iRe(\xi \overline{z})} f( \xi,w)d \xi_1d\xi_2$$ the inverse Fourier transform.

\begin{lemma}\label{Lemma 3.2}     Putting $g = g(z,w) = {\mathcal F}_z^{-1} (f )(z,w)$ we have:
\begin{enumerate}
\item
$$\partial_z g = \frac{i}{2}\overline\xi g \    ;  \partial_z^{r} g = {(\frac{i}{2}\overline\xi)}^r g, r=2,3,\dots $$
\item
 $$ \partial_{\overline z} g = \frac{i}{2}\xi g \    ;  \partial_{\overline z}^{r} g = {(\frac{i}{2}\xi)}^r g, r=2,3,\dots $$
\item 
$${\mathcal F}_z(zg) = 2i\partial_{\overline\xi}{\mathcal F}_z(g) = 2i\partial_{\overline \xi}f \    ;  {\mathcal F}_z(\overline{z}g) = 2i\partial_{\xi}{\mathcal F}_z(g) = 2i\partial_{\xi}f $$    
\item
$$ \partial_w g = \partial_w ({\mathcal F}_z^{-1}(f)) = {{\mathcal F}_z} ^{-1}(\partial_{w}f);\     \partial_{\overline w}g = \partial_{\overline w} ({\mathcal F}_z^{-1}(f)  = {{\mathcal F}_z}^{-1}(\partial_{\overline w}f)$$
\end{enumerate}  
\end{lemma}

We also need another Lemma which will be used in the sequel.

\begin{lemma}\label{Lemma 3.3}     With $g = {\mathcal F}_z^{-1}$$(f)($$z,w)$, we have: 
\begin{enumerate}
\item
$$    {\mathcal F}_z(P^0(\widetilde{A},g)) = i(\alpha \partial_{\overline \xi} + \overline {\alpha} \partial_{\xi})f + \frac{1}{ 2} \beta e^w f + \frac{1}{2} \overline {\beta} e^{\overline w} f. $$
\item
$$   {\mathcal F}_z(P^1(\widetilde{A},g)) = \overline {\alpha} \partial_{\overline w}f + \alpha \partial_{w}f - \overline {\beta} e^{\overline w} (\frac{i}{2} \xi)f - \beta e^w (\frac{i}{2}\overline {\xi})f. $$
\item
$$    {\mathcal F}_z(P^r(\widetilde{A},g)) = {(-1)}^r.2^{r-1}[\overline \beta{ e^{\overline w}} (\frac{i}{2}\xi)^r + \beta e^w (\frac{i}{2}\overline \xi)^r]f \ \ \ \ \ \      \forall r \ge 2. $$
\end{enumerate}
\end{lemma}

\begin{proposition}\label{Proposition 3.4}
For each $A = \left(\begin{matrix}\alpha & \beta \cr 0 & 0 \cr\end{matrix}\right) \in \aff({\bf C}) $
 and for each compactly supported $C^{\infty}$-function $f \in C_0^{\infty}({\bf C} \times {\bf H}_k)$, we have:
\begin{equation} {\ell}_A^{(k)}{f} := {\mathcal F}_z \circ \ell_A^{(k)} \circ {\mathcal F}_z^{-1}(f) = [\alpha (\frac{1}{2} \partial_w - \partial_{\overline \xi})f + \overline \alpha (\frac{1 }{2}\partial_{\overline w} - \partial_\xi)f + \end{equation}
$$+\frac{i}{2}(\beta e^{w-\frac{1}{2}\overline \xi} + \overline \beta e^{\overline w - \frac{1}{2} \xi})f] $$
\end{proposition}

\begin{remark}\label{Remark 3.5}{\rm Setting new variables  u = $w - \frac{1}{ 2}\overline{\xi}$;$v = w + \frac{1 }{2}{\overline{\xi}}$ we have
\begin{equation}\hat {\ell}_A^{(k)}(f) = \alpha\frac{ \partial f }{\partial u}+ \overline{\alpha}\frac{\partial f }{\partial{\overline{u}}}+ \frac{i }{2}(\beta e^{u}+\overline{\beta}e^{\overline{u}})f \vert_{(u,v)}\end{equation}
i.e $\hat {\ell}_A^{(k)} = \alpha\frac{ \partial }{\partial u}+ \overline{\alpha}\frac{ \partial  }{\partial{\overline{u}}}+ \frac{i }{2}(\beta e^{u}+\overline{\beta}e^{\overline{u}})$,which provides a ( local) representation of the Lie algebra  aff({\bf C}).
}\end{remark}

{\it The Irreducible Representations of $\widetilde{\Aff}({\bf C})$. }
Since $\hat {\ell}_A^{(k)}$ is a representation of the Lie algebra  $\widetilde{\hbox {Aff}} ({\bf C})$, we have:
$$\exp(\hat {\ell}_A^{(k)}) = \exp\bigl(\alpha\frac{ \partial }{\partial {u}}+ \overline{\alpha}\frac{ \partial  }{\partial{\overline{u}}}+ \frac{i }{2}(\beta e^{u}+\overline{\beta}e^{\overline{u}})\bigr)$$ is just the corresponding representation of the corresponding connected and simply connected Lie group $\widetilde\Aff ({\bf C})$.

Let us first recall the well-known list of all the irreducible unitary representations of the group of affine transformation of the complex straight line, see [D] for more details.

\begin{theorem}\label{Theorem 4.1}
Up to unitary equivalence, every irreducible unitary representation of $\widetilde{\hbox {Aff}} ({\bf C})$ is unitarily equivalent to one the following one-to-another nonequivalent irreducible unitary representations:
\begin{enumerate}
\item The unitary characters of the group, i.e the one dimensional unitary representation $U_{\lambda},\lambda \in {\bf C}$, acting in ${\bf C}$ following the formula
$U_{\lambda}(z,w) = e^{{i\Re(z\overline{\lambda})}}, \forall (z,w) \in \widetilde{\Aff} ({\bf C}), \lambda \in {\bf C}.$
\item The infinite dimensional irreducible representations $T_{\theta},\theta \in {\mathbf S}^1$, acting on the Hilbert space $L^{2}(\mathbf R\times \mathbf S ^1)$ following the formula:
\begin{equation}\Bigr[T_{\theta}(z,w)f\Bigl](x) = \exp \Bigr(i(\Re(wx)+2\pi\theta[\frac{\Im(x+z) }{2\pi}])\Bigl)f(x\oplus z),\end{equation}
Where \ $(z,w) \in\widetilde{\Aff}({\bf C})$  ;  $x \in {\bf R}\times {\mathbf S} ^1= {\bf C} \backslash \{0\}; f \in L^{2}({\bf R}\times {\mathbf S} ^1);$
$$ x\oplus z = Re(x+z) +2 \pi i \{\frac{\Im(x+z) }{2\pi}\}$$
\end{enumerate}
\end{theorem}
In this section we will prove the following important Theorem which
is very interesting for us both in theory and practice.
\par
\begin{theorem}\label{Theorem 4.2}
The representation $\exp(\hat {\ell}_A^{(k)})$ of the group $\widetilde{\Aff}({\bf C})$ is the irreducible unitary representation 
$T_\theta$ of $\widetilde{\Aff}({\bf C})$ associated, following the orbit method construction, to the orbit $\Omega$, i.e:
$$\exp(\hat {\ell}_A^{(k)})f(x) = [T_\theta (\exp A)f](x),$$
where $f \in L^{2}({\bf R}\times {\mathbf S} ^1) ; A = \begin{pmatrix}\alpha & \beta \cr 0 & 0 \cr\end{pmatrix} \in \aff({\bf C}) ; \theta \in {\mathbf S}^1 ; k = 0, \pm1,\dots$
\end{theorem} 

\begin{remark}\label{Remark 4.3} {\rm
We say that a real Lie algebra ${\mathfrak g}$ is in the class $\overline{MD}$ if every K-orbit is of dimension, equal 0 or dim ${\mathfrak g}$. Further more, one proved that
([D, Theorem 4.4]) 
Up to isomorphism, every Lie algebra of class $\overline {MD}$ is one of the following:
\begin{enumerate}
\item Commutative Lie algebras.
\item Lie algebra $\aff({\bf R})$ of affine transformations of the real straight line
\item Lie algebra $\aff({\bf C})$ of affine transformations of the complex straight line.
\end{enumerate}
Thus, by calculation for the group of affine transformations of the real straight line $\Aff({\bf R})$ in [DH] and here for the group affine transformations of the complex straight line $\Aff({\bf C})$ we obtained  a description of the quantum $\overline {MD}$ co-adjoint orbits.
}\end{remark}

\subsection{$MD_4$-groups}
We refer the reader to the results of Nguyen Viet Hai \cite{hai3}-\cite{hai4} for the class of $MD_4$-groups (i.e. 4-dimensional solvable Lie groups, all the coadjoint of which are of dimension 0 or maximal). It is interesting that here he obtained the same exact computation for $\star$-products and all representations.

\subsection{$SO(3)$} As an typical example of compact Lie group, the author proosed Job A. Nable to consider the case of $SO(3)$. We refer the reader to the  results of Job Nable \cite{nable1}-\cite{nable3}. In these examples, it is interesting that the $\star$-products, in some how as explained in these papers, involved the Maslov indices and Monodromy Theorem.

	\subsection{Exponetial groups}
Arnal-Cortet constructed star-products for this case \cite{arnalcortet1}-\cite{arnalcortet2}.
	\subsection{Compact groups}
We refer readers to the works of C. Moreno \cite{moreno}.

\section{Algebraic Noncommutative Chern Characters}

\subsection{Modification for inductive limits}

Let $G$ be a compact group, $\HP_*(C^*(G))$ the periodic cyclic homology
introduced
in \S2. Since $C^*(G) = \varinjlim_N \prod_{i=1}^N \Mat_{n_i}({\mathbf C})$,
$\HP_*(C^*(G))$ coincides with the $\HP_*(C^*(G))$ defined by J. Cuntz-D.
Quillen
\cite{CQ}.

\begin{lem}\label{31}
Let $\{I_N\}_{N\in {\mathbf N}}$ be the above defined collection of ideals
in $C^*(G)$. Then $$K_*(C^*(G)) = \varinjlim_{N\in {\mathbf N}} K_*(I_N) = K_*({\mathbf C}(\mathbf T)),$$ where $\mathbf T$ is the fixed maximal torus in $G$.
\end{lem}
 First note that the algebraic K-theory of C*-algebras has the stability property $$K_*(A \otimes M_n(\mathbf C)) \cong K_*(\mathbf C(\mathbf T)).$$ Hence, $$\varinjlim K_*(I_{n_i}) \cong K_*(\prod_{w=\mbox{ highest weight}} {\mathbf C}_w)\cong K_*(\mathbf C(\mathbf T)),$$ by Pontryagin duality.

\subsection{An algebraic construction of noncommutative Chern characters}

J. Cuntz and D. Quillen \cite{CQ} defined the so called $X$-complexes of
${\mathbf C}$-algebras and then used some ideas of Fedosov product to
define algebraic Chern
characters. We now briefly recall the their definitions. For a
(non-commutative) associate ${\mathbf C}$-algebra $A$, consider the space
of even
non-commutative differential forms $\Omega^+(A) \cong RA$, equipped with
the Fedosov
product
$$\omega_1 \circ \omega_2 := \omega_1\omega_2 - (-1)^{\vert \omega_1\vert}
d\omega_1
d\omega_2,$$ see \cite{CQ}. Consider also the ideal $IA := \oplus_{k\geq 1}
\Omega^{2k}(A)$. It is easy to see that $RA/IA \cong A$ and that $RA$
admits the
universal property that any based linear map $\rho : A \to M$ can be
uniquely extended to a
derivation $D : RA \to M$. The derivations $D : RA \to M$ bijectively
correspond to lifting
homomorpisms from $RA$ to the semi-direct product $RA \oplus M$, which also
bijectively correspond to linear map $\bar\rho : \bar{A}= A /{\mathbf C}
\to M$ given by $$
a\in \bar{A} \mapsto D(\rho a).$$ From the universal property of
$\Omega^1(RA)$, we
obtain a bimodule isomorphism $$RA \otimes \bar{A} \otimes RA \cong
\Omega^1(RA).$$ As in \cite{CQ}, let $\Omega^-A = \oplus_{k\geq 0}
\Omega^{2k+1}A$. Then we have $$\Omega^{-}A \cong RA \otimes \bar{A} \cong
\Omega^1(RA)_\# := \Omega^1(RA)/[(\Omega^1(RA),RA)].$$
\par
J. Cuntz and D. Quillen proved
\begin{thm}(\cite{CQ}, Theorem1):
There exists an isomorphism of
${\mathbf Z}/(2)$-graded complexes
$$\Phi : \Omega A = \Omega^+A \oplus \Omega^{-}A \cong RA \oplus
\Omega^1(RA)_\#,$$ such that
$$\Phi : \Omega^+A \cong RA,$$ is defined by $$\Phi(a_0da_1\dots da_{2n} =
\rho(a_1)\omega(a_1,a_2) \dots \omega(a_{2n-1},a_{2n}),$$
and $$ \Phi : \Omega^{-}A \cong \Omega^1(RA)_\#,$$ $$\Phi(a_0da_1\dots
da_{2n+1})
= \rho(a_1)\omega(a_1,a_2)\dots \omega(a_{2n-1},a_{2n})\delta(a_{2n+1}).$$
With
respect to this identification, the product in $RA$ is just the Fedosov
product on even
differential forms and the differentials on the $X$-complex
$$X(RA) : \qquad RA\cong \Omega^+A \to \Omega^1(RA)_\# \cong \Omega^{-}A
\to RA
$$ become the operators
$$\beta = b - (1+\kappa)d : \Omega^{-}A \to \Omega^+A,$$ $$\delta =
-N_{\kappa^2} b
+ B : \Omega^+A \to \Omega^{-}A,$$ where $N_{\kappa^2} =
\sum_{j=0}^{n-1} \kappa^{2j}$, $\kappa(da_1\dots da_n) := da_n\dots da_1$.
\end{thm}
\par
Let us denote by $IA \triangleleft RA$ the ideal of even non-commutative
differential forms
of
order $\geq 2$. By the universal property of $\Omega^1$ $$\Omega^1(RA/IA) =
\Omega^1RA/((IA)\Omega^1RA + \Omega^1RA.(IA) + dIA).$$ Since $\Omega^1RA =
(RA)dRA = dRA.(RA)$, then $\Omega^1RA(IA) \cong IA\Omega^1RA \mod
[RA,\Omega^1R].$
$$\Omega^1(RA/IA)_\# = \Omega^1RA /([RA,\Omega^1RA]+IA.dRA + dIA).$$ For
$IA$-adic tower $RA/(IA)^{n+1}$, we have the complex $${\mathcal
X}(RA/(IA)^{n+1}) : \qquad RA/IA^{n+1} \leftarrow
\Omega^1RA/([RA,\Omega^1RA]+(IA)^{n+1}dRA + d(IA)^{n+1}).$$
Define
$${\mathcal X}^{2n+1}(RA,IA) : \quad RA/(IA)^{n+1} \to
\Omega^1RA/([RA,\Omega^1RA]+(IA)^{n+1}dRA + d(IA)^{n+1}) $$ $$\to
RA/(IA)^{n+1},$$
$${\mathcal X}^{2n}(RA,IA): \quad RA/((IA)^{n+1} +[RA,IA^n]) \to
\Omega^1RA/([RA,\Omega^1RA]+d(IA)^ndRA)$$ $$\to RA/((IA)^{n+1}
+[RA,IA^n]).$$
One has
$$b((IA)^ndIA) = [(IA)^n,IA] \subset (IA)^{n+1},$$ $$d(IA)^{n+1} \subset
\sum_{j=0}^n (IA)^jd(IA)(IA)^{n-j} \subset (IA)^n dIA
+ [RA,\Omega^1RA].$$
and hence
$${\mathcal X}^1(RA,IA = X(RA,IA),$$
$${\mathcal X}^0(RA,IA) = (RA/IA)_\#.$$
There is an sequence of maps between complexes $$\dots \to X(RA/IA) \to
{\mathcal
X}^{2n+1}(RA,IA)\to {\mathcal X}^{2n}(RA,IA) \to X(RA/IA) \to \dots $$ We
have the
inverse limits
$$\hat{X}(RA,IA) := \varprojlim X(RA/(IA)^{n+1}) = \varprojlim {\mathcal
X}^n(RA,IA).$$
Remark that
$${\mathcal X}^q = \Omega A/F^q\Omega A,$$ $$\hat{X}(RA/IA)= \hat{\Omega}A.$$

We quote the second main result of J. Cuntz and D. Quillen (\cite{CQ},
Thm2), namely:

$$H_i\hat{\mathcal X}(RA,IA) = \HP_i(A).$$

We now apply  this machinery to our case. First we have the following.
\begin{lem}
$$\varinjlim \HP^*(I_N) \cong \HP^*({\mathbf C}({\mathbf T})).$$ \end{lem}
 By similar arguments as in the previous lemma \ref{31}. More
precisely, we
have $$\HP(I_{n_i}) = \HP(\prod_{w=\mbox{heighest weight}} {\mathbf C}_w)
\cong
\HP({\mathbf C}({\mathbf T}))$$ by Pontryagin duality.

Now, for each idempotent $e\in M_n(A)$ there is an unique element $x\in
M_n(\widehat{RA})$.
Then the element $$\tilde{e} := x + (x-\frac{1}{2})\sum_{n\geq 1} \frac{2^n(2n-
1)!!}{n!}(x-x^{2n})^{2n}\in M_n(\widehat{RA})$$ is a lifting of $e$ to an
idempotent
matrix in $M_n(\widehat{RA})$. Then the map $[e] \mapsto tr(\tilde{e})$
defines the map
$K_0 \to H_0(X(\widehat{RA})) = \HP_0(A)$. To an element $g\in\GL_n(A)$ one
associates an element $p\in \GL(\widehat{RA})$ and to the element $g^{-1}$
an element
$q\in \GL_n(\widehat{RA})$ then put
$$x = 1- qp, \mbox{ and } y = 1-pq.$$
And finally, to each class
$[g]\in \GL_n(A)$ one associates $$tr(g^{-1}dg) = tr(1-x)^{-1}d(1-x) =
d(tr(log(1-x))) =
-tr\sum_{n=0}^\infty x^ndx\in \Omega^1(A)_\#.$$ Then $[g] \to tr(g^{-1}dg)$
defines
the map $K_1(A) \to HH_1(A) = H_1(X(\widehat{RA})) = \HP_1(A)$.

\begin{defn} Let $\HP(I_{n_i})$ be the periodic cyclic cohomology defined by
Cuntz-Quillen. Then the pairing
$$K_*^{alg}(C^*(G)) \times \bigcup_N \HP^*(I_N) \to {\mathbf C}$$ defines an
algebraic non-commutative Chern character $$ch_{alg} : K_*^{alg}(C^*(G)) \to
\HP_*(C^*(G)),$$ which gives us a variant of non-commutative Chern
characters with
values in $\HP$-groups. \end{defn}

We close this section with an algebraic analogue of theorem 3.2.

\begin{thm}
Let $G$ be a compact group and ${\mathbf T}$ a fixed maximal compact torus
of $G$.
Then in the notations of 4.3, the Chern character $$ch_{alg} : K_*(C^*(G)) \to
\HP_*(C^*(G))$$ is an isomorphism, which can be identified with the
classical Chern
character $$ch: K_*({\mathbf C}({\mathbf T})) \to \HP({\mathbf C}({\mathbf
T}))$$
which is also an isomorphism.
\end{thm}

\section*{Acknowledgments}

This survey is an extended version of  invited lectures at the {\it Operator Theory Seminar}, the University of Iowa, September 14 and 21, 1999 and at the {\it INFAS (Iowa-Nebraska Funtional Analysis Seminar)}, Drake University, October 30, 1999 and was completed during the stay of the author as a visiting
mathematician at the Department of mathematics, The University of Iowa. The author would like to express the deep and sincere thanks to Professor Tuong Ton-That and his spouse, Dr. Thai-Binh Ton-That for their effective helps and kind attention they provided during the stay in Iowa, and also for a discussion about the PBW Theorem. The deep thanks are also addressed to the organizers of the Seminar on Mathematical Physics, Seminar on Operator Theory in Iowa and the Iowa-Nebraska Functional Analysis Seminar (INFAS), in particular the professors Raul Curto, Palle Jorgensen, Paul Muhly and Tuong Ton-That for the stimulating scientific atmosphere.
 
The author would like to thank the University of Iowa for the hospitality and the scientific support, the Alexander von Humboldt Foundation, Germany, for an effective support.

{\noindent\sc Institute of Mathematics,  National Centre for Science and Technology}, {\rm  P. O. Box 631, Bo Ho, VN-10000, Hanoi, Vietnam}\\ 
Email: {\tt dndiep@hn.vnn.vn} 

{\sc{\rm and}\\ Department of Mathematics, The University of Iowa}, {\rm 14 MacLean Hall, Iowa City, IA 52242-1419, U.S.A.}\\
Email: {\tt ndiep@math.uiowa.edu}

\end{document}